\documentclass[a4paper,10pt,leqno,twoside]{tobart}

\usepackage[english]{babel} \usepackage{inputenc, amsmath, amssymb, latexsym,
  epic, epsfig, rotating, fancyheadings, amsthm, pifont, empheq}

\newcommand{\follows}{\ensuremath{\Rightarrow}}

\newcommand{\ld}{\ensuremath{,\ldots,}}
\newcommand{\ssq}{\ensuremath{\subseteq}}
\newcommand{\smin}{\ensuremath{\setminus}}
\newcommand{\eps}{\ensuremath{\varepsilon}}

\newcommand{\T}{\ensuremath{\mathbb{T}}}

\newcommand{\N}{\ensuremath{\mathbb{N}}} 
\newcommand{\R}{\ensuremath{\mathbb{R}}}
\newcommand{\Z}{\ensuremath{\mathbb{Z}}}
\newcommand{\Q}{\ensuremath{\mathbb{Q}}}

\newcommand{\supp}{\ensuremath{\mathrm{supp}}}

\newcommand{\inte}{\ensuremath{\mathrm{int}}}

\newcommand{\diam}{\ensuremath{\mathrm{diam}}}
\newcommand{\spann}{\ensuremath{\mathrm{span}}}
\newcommand{\conv}{\ensuremath{\mathrm{Conv}}}

\newcommand{\ntorus}[1][2]{\ensuremath{\mathbb{T}^{#1}}}

\newcommand{\torus}{\ensuremath{\mathbb{T}^2}}
\newcommand{\homd}{\ensuremath{\mathrm{Homeo}_0(\mathbb{T}^d)}}
\newcommand{\homeo}{\ensuremath{\mathrm{Homeo}}}
\newcommand{\homtwo}{\ensuremath{\mathrm{Homeo}_0(\mathbb{T}^2)}}

\newcommand{\alphlist}{\begin{list}{(\alph{enumi})}{\usecounter{enumi}\setlength{\parsep}{2pt}
      \setlength{\itemsep}{1pt} \setlength{\topsep}{3pt}
      \setlength{\partopsep}{2pt}}}

\newcommand{\arablist}{\begin{list}{(\arabic{enumi})}{\usecounter{enumi}\setlength{\parsep}{2pt}
          \setlength{\itemsep}{1pt} \setlength{\topsep}{3pt}
          \setlength{\partopsep}{2pt}}}

\newcommand{\romanlist}{\begin{list}{(\roman{enumi})}{\usecounter{enumi}\setlength{\parsep}{2pt}
              \setlength{\itemsep}{1pt} \setlength{\topsep}{3pt}
              \setlength{\partopsep}{2pt}}}

 \newcommand{\listend}{\end{list}}

\newcommand{\bulletlist}{\begin{list}{$\bullet$}{\setlength{\parsep}{2pt}
                \setlength{\itemsep}{1pt} \setlength{\topsep}{3pt}
                \setlength{\partopsep}{2pt}\setlength{\leftmargin}{15pt}}}

\newcommand{\roundqed}{{\hfill \Large $\circ$}}

\newcommand{\myproof}{\textit{Proof. }}

\newcommand{\foot}{\footnote}

\newcommand{\kcup}{\ensuremath{\bigcup_{k\in\N}}}

\newcommand{\nLim}{\ensuremath{\lim_{n\rightarrow\infty}}}
\newcommand{\iLim}{\ensuremath{\lim_{i\rightarrow\infty}}}

\newcommand{\kLim}{\ensuremath{\lim_{k\rightarrow\infty}}}

\newcommand{\ntel}{\ensuremath{\frac{1}{n}}}

\newcommand{\halb}{\ensuremath{\frac{1}{2}}}

\newcommand{\viertel}{\ensuremath{\frac{1}{4}}}

\newcommand{\dreiviertel}{\ensuremath{\frac{3}{4}}}

\newtheoremstyle{tobthm}{3pt}{3pt}{\itshape}{0pt}{\bfseries}{.}{0.5eM}{}
\theoremstyle{tobthm}
\newtheorem{definition}{Definition}[section]
\newtheorem{thm}[definition]{Theorem}

\newtheorem{lem}[definition]{Lemma}
\newtheorem{cor}[definition]{Corollary}  
\newtheorem{prop}[definition]{Proposition}

\newtheorem{claim}[definition]{Claim}

\newtheoremstyle{tobrem}{3pt}{3pt}{\normalfont}{0pt}{\bfseries}{.}{0.5em}{}
\theoremstyle{tobrem}
\newtheorem{rem}[definition]{Remark}

\numberwithin{equation}{section}

\newcommand{\vlin}{\ensuremath{\lambda v + \{v\}^\perp}}

\newcommand{\Llv}{\ensuremath{L_{\lambda,v}}}
\newcommand{\vbr}[1]{\ensuremath{\left\langle #1 \right\rangle_v}}
\newcommand{\vpbr}[1]{\ensuremath{\left\langle #1 \right\rangle_{v^\perp}}}
\newcommand{\wh}{\ensuremath{\widehat}}

\title{\Large\textsc{Elliptic stars in a chaotic night}}
\author{T.~J\"ager\thanks{TU Dresden, Institut f\"ur Analysis. Email:
    {\tt Tobias.Oertel-Jaeger@tu-dresden.de}}}

\begin{document}

\setlength{\abovedisplayskip}{0.8ex}
\setlength{\abovedisplayshortskip}{0.6ex}

\setlength{\belowdisplayskip}{0.8ex}
\setlength{\belowdisplayshortskip}{0.6ex}

\maketitle 

\abstract{A recurrent theme in the description of phase portraits of dynamical
  systems is that of {\em elliptic islands in a chaotic sea}. Usually this
  picture is invoked in the context of smooth twist maps of the annulus or the
  torus, like the standard map. In this setting `elliptic islands' refers to the
  topological disks bounded by periodic smooth curves surrounding elliptic
  periodic points. Establishing the existence of these curves is one of the many
  achievements of KAM-theory.

  The aim of this note is to approach the topic from a different
  angle, namely from the viewpoint of rotation theory in a purely
  topological setting.  We study homeomorphisms of the two-torus,
  homotopic to the identity, which have no wandering open sets (as in
  the area-preserving case) and whose rotation set has non-empty
  interior. We define local rotation subsets $\rho_F(U)$ by
  restricting Misiurewicz and Ziemian's definition of the rotation set
  to starting points in a small open disk $U$. Our main result is the
  following dichotomy: Either $\rho_U(F)$ is reduced to a single
  rational vector and $U$ is contained in a periodic topological open
  disk which contains a periodic point, or $\rho_U(F)$ is large, in
  the sense that its convex hull has non-empty interior. This allows
  to distinguish an `elliptic' and a `chaotic' regime, and as a
  consequence we obtain that in the chaotic region the dynamics
  are sensitive with respect to initial conditions.

  In order to demonstrate these results we introduce a parameter family of
  smooth toral diffeomorphisms that is inspired by an example of Misiurewicz and
  Ziemian. The pictures obtained from simulations in this family motivate an
  alternative formulation of the original theme. }

\enlargethispage*{1000pt}

\begin{figure}[h] 
\begin{center}
\epsfig{file=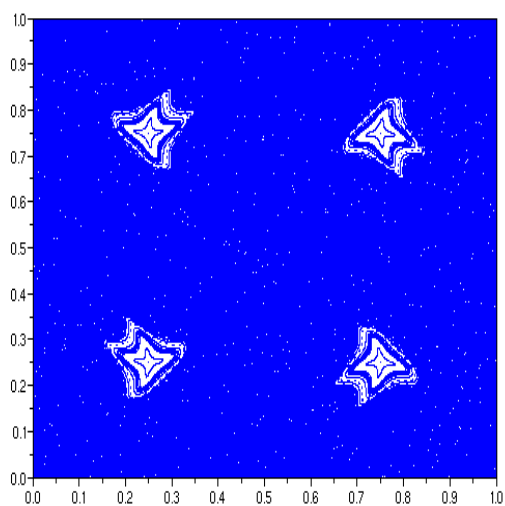, width=0.58\linewidth, height=0.58\linewidth}
\end{center} {\caption{ \small Elliptic islands surrounding two 2-periodic
    orbits of the map $f_\alpha(x,y) =
    (x+\alpha\sin(2\pi(y+\alpha\sin(2\pi x))),y+\alpha\sin(2\pi x))$ with
    $\alpha=0.5$.} \label{f.1}}
\end{figure}

\pagebreak

\section{Introduction}

We denote by \homtwo\ the set of homeomorphisms of the two-torus that are
homotopic to the identity. Given a lift $F:\R^2 \to \R^2$ of $f\in\homtwo$,
Misiurewicz and Ziemian \cite{misiurewicz/ziemian:1989} introduced the rotation
set of $F$ as
\begin{equation}\label{e.rot-set}
  \rho(F) \ := \ \left\{ \rho \in \R^2 \left|\  \exists n_i\nearrow\infty,\ z_i\in\R^2 :
 \iLim \left(F^{n_i}\left(z_i\right) - z_i\right)/n_i \ = \ \rho \right.\right\} \ .
\end{equation}
This set is always compact and convex \cite{misiurewicz/ziemian:1989}.
Further, the properties of $\rho(F)$ have strong implications for the
dynamics of $f$. In particular, this is true for the situation we will
concentrate, namely when $\rho(F)$ has non-empty interior. In this
case all rotation vectors in $\inte(\rho(F))$ are realised on minimal
sets \cite{franks:1989,misiurewicz/ziemian:1991,jaeger:2009b} and the
topological entropy of $f$ is strictly positive
\cite{llibre/mackay:1991}. Further the set $\left\{f\in\homtwo\mid
  \inte(F)\neq \emptyset\right\}$ is an open and therefore, in a
topological sense, large subset of \homtwo\
\cite{misiurewicz/ziemian:1991}.

Our aim is to give some meaning to the notion of {\em elliptic islands in a
  chaotic sea} in this purely topological setting. To that end, we restrict the
definition in (\ref{e.rot-set}) to orbits starting in some subset
$U\ssq\torus$. Let $\pi:\R^2\to\T^2$ denote the canonical projection. We define
the {\em rotation subset on $U$} by
\begin{equation} \label{e.rotation-subset} \rho_U(F) \ := \ \left\{ \rho \in
    \R^2 \left| \ \exists z_i\in \pi^{-1}(U),\ n_i\nearrow \infty : \iLim
      \left(F^{n_i}\left(z_i\right)-z_i\right)/n_i = \rho \right. \right\} \ .
\end{equation}
In general, even when $U$ is open $\rho_U(F)$ can be much smaller than
$\rho(F)$. For instance, when $f$ is a sufficiently smooth toral diffeomorphism
then generic elliptic periodic points are surrounded by periodic invariant
curves (see, for example, \cite{moser:1962,katok/hasselblatt:1997}). The
rotation subsets of the corresponding topological disks contain a single
rational rotation vector, whereas $\rho(F)$ may have non-empty interior. A more
general example is sketched in Remark~\ref{r.denjoy-construction} below.

However, when $U$ is open and recurrent, then in a number of situations
$\rho(F)$ is already determined by $\rho_U(F)$. In order to give precise
statements, we need some notation and terminology. We say $U\ssq \torus$ is {\em
  bounded} if the connected components of its lift to $\R^2$ are bounded. Given
$f\in\homtwo$ we say $U$ is {\em wandering} if $f^n(U) \cap U = \emptyset \
\forall n\geq 1$ and {\em non-wandering} otherwise. We call $U$ {\em recurrent}
if there exist infinitely many $n\in\N$ with $f^n(U) \cap U \neq \emptyset$.  We
call $z\in\torus$ {\em wandering} if it is contained in some wandering open set
and {\em non-wandering} otherwise. It is easy to see that if $U$ is open and
contains a non-wandering point then it is recurrent. Finally, we say that $f$
is {\em non-wandering} if it has all points non-wandering. In this case all open
sets are recurrent. Note that any area-preserving toral homeomorphism is
non-wandering. 

Given a lift $F$ of $f\in\homtwo$, let $\varphi_n(z) = \left(F^n(z)-z\right)/n$. If
$\lambda\in\R$ and $v\in\R^2\smin\{0\}$, let $\Llv = \vlin$.  
\begin{thm} \label{t.semilocal-rotsets} Suppose $F$ is a lift of $f\in\homtwo$
  and $U\ssq \ntorus$ is open, bounded, connected and recurrent.  \alphlist
\item If $\rho_U(F) \ssq \Llv$ for some $\lambda\in\R$, $v\in\R^2\smin\{0\}$, then
  either $\rho(F) \ssq \Llv$ or $\rho_U(F)$ is reduced to a single rational
  vector.
\item If ${\cal S}$ is a line segment of positive length without rational points
  and $\rho_U(F) = {\cal S}$, then $\rho(F) = {\cal S}$. Further, $\varphi_n(U)$
  converges to ${\cal S}$ in Hausdorff distance as $n\to \infty$.
\item If $\rho_U(F) = \{\rho\}$ with $\rho\in\R^2$ irrational\foot{We call
   a vector $\rho=\left(\rho_1,\rho_2\right)\in \R^2$ {\em rational} if $\rho\in\Q^2$, {\em
      irrational} if $\rho_1,\rho_2,\rho_1/\rho_2 \notin\Q$ and {\em
      semi-rational} if it is neither rational nor irrational.} then
  $\rho(F)=\{\rho\}$.
\listend
 \end{thm}

 \begin{rem} \label{r.denjoy-construction} We note that when no recurrence
   assumption is made no relation between $\rho_U(F)$ and $\rho(F)$ can be expected.
   Without going into detail, we want to mention a possible way to construct
   respective examples: When $\rho(F)$ has non-empty interior, then for any
   compact connected subset $C \ssq \rho(F)$ there exists a point $z\in\T^2$
   with $\rho_{\{z\}}(F) = C$ \cite{llibre/mackay:1991}. By blowing up the
   points in the orbit of $z$ to small disks in a Denjoy-like construction one
   may thus obtain a wandering open set $U$ whose rotation set is an arbitrary
   compact connected subset of the rotation set.
\end{rem}

Theorem~\ref{e.rot-set}(i) implies that if $f$ is non-wandering and
$\rho(F)$ has non-empty interior then the rotation subset of $U$ can
only be contained in a line if it is reduced to a single rational
rotation vector, that is, $\rho_U(F) = \{\rho\}$ with $\rho\in\Q^2$.
Together with some additional details on the rational case, this
yields our main result.  \enlargethispage*{1000pt}
 \begin{thm}
   \label{t.elliptic-islands} Let $F:\R^2 \to \R^2$ be a lift of $f\in\homtwo$
   and suppose that $f$ is non-wandering and $\rho(F)$ has non-empty interior.
   Then for any open, bounded and connected set $U$ one of the following two holds.  \romanlist
 \item $\rho_U(F)$ is reduced to a single rational vector $\rho$ and $U$ is
   contained in an embedded topological open disk $D\ssq \torus$ which is
   invariant under some iterate $f^p$ and contains a $p$-periodic point.
 \item The convex hull of $\rho_U(F)$ has non-empty interior.
   \listend
 \end{thm}
\pagebreak

The above result allows to give an intrinsic definition of `elliptic' and
`chaotic' regions. Given a set $A\ssq\R^2$ we denote by $\conv(A)$ its convex hull and
by $\inte(A)$ its interior.
\begin{definition}
  Suppose $f\in\homtwo$. Let
  \begin{eqnarray*} {\cal E}(f) & := & \left\{ z\in\torus \mid \#\rho_U(F) = 1 \
    \textrm{ for some open
      neighbourhood } U \textrm{ of } x \right\} \quad , \\
    {\cal C}(f) & := & \left\{ z\in\torus \mid \inte\left(\conv(\rho_U(F))\right) \neq
    \emptyset \ \ \forall \textrm{open neighbourhoods } U \textrm{ of } x \right\} \ .
  \end{eqnarray*}
\end{definition}
A point $z\in\torus$ is called {\em $\eps$-Lyapunov stable} if there exists some
$\delta > 0$ such that $f^n\left(B_\delta(z)\right) \ssq
B_\eps\left(f^n(z)\right) \ \forall n\in\N$ and {\em Lyapunov stable} if it is
$\eps$-Lyapunov stable for all $\eps > 0$. As one should expect for a notion of
stability, Lyapunov stable points do not occur in the `chaotic' regime.
\begin{prop} \label{p.chaotic} Suppose $f\in\homtwo$ and $\rho(F)$ has non-empty
  interior.\alphlist
\item If $f$ is non-wandering then no point in ${\cal C}(f)$ is
  $\frac{1}{2}$-Lyapunov stable.
\item If $f$ is area-preserving, $U$ is a connected and bounded neighbourhood of
  $z\in {\cal C}(f)$ and $\widehat U$ is a connected component of $\pi^{-1}(U)$
  then $\limsup_{n\to\infty}\diam\ntel\left(F^n(\widehat U)\right) > 0$. \listend
\end{prop}
Note that in contrast to this, in the construction sketched in
Remark~\ref{r.denjoy-construction} all points in the wandering topological disks
will be Lyapunov stable provided the diameter of these disks goes to zero along
the orbit.  \medskip

The paper is organised as follows: In Section~\ref{Basic} we collect a
number of basic statements on rotation subsets.
Section~\ref{Semilocal-rotsets} then contains the technical core of
the paper. We work on the universal cover $\R^2$ and consider bounded
open and connected sets that intersect their image. In this setting,
we describe a number of situations in which the rotation subset
already determines the rotation set, or at least forces it to be
contained in a line segment. In Section~\ref{Proofs} these statements
are then used to prove the main results. Finally, in
Section~\ref{MZ-family} we provide some explicit examples to which our
results apply. To that end, we introduce a parameter family of smooth
torus diffeomorphisms that is based on an example by Misiurewicz and
Ziemian in \cite{misiurewicz/ziemian:1991}. For appropriate parameter
values these maps have a rotation set with non-empty interior, and
simulations clearly indicate that they also exhibit elliptic islands.

\medskip

{\bf Acknowledgements.} This work was supported by an Emmy Noether Grant (Ja
1721/2-1) of the German Research Council (DFG). The results were partially
presented at the Workshop on Dynamics in Dimension two, April 2009 in Puc\'on
(Chile). I would like to thank the organisers Mario Ponce and Andres Navas for
the great opportunity they created. I am further indebted to Andres Koropecki
and Patrice Le Calvez for thoughtful remarks and stimulating discussions and to
Hendrik Vogt for drawing my attention to the fine structure of the elliptic
islands depicted in Figure~\ref{fig.7}. Finally, I would like to thank
Jean-Christophe Yoccoz for his support during my time at the Coll\`ege de
France, where a part of this work was carried out.

\section{Some basic results on rotation subsets} \label{Basic}

The aim of this section is to collect a number of elementary
statements on rotation subsets and rotation vectors that will be used
in the later sections. For the purposes of this section there is no
need to restrict to dimension 2. Hence, we will work on $\T^d$
($d\in\N$), with the definitions of the rotation set and rotation
subsets analogous to those on $\torus$. 
\smallskip 

\noindent {\em Notation.} \ We denote the Euclidean scalar product of
vectors $v,w\in\R^d$ by $\langle v,w\rangle$ and also write $\vbr{w}$
instead of $\langle v,w \rangle$. By $\|v\|=\sqrt{\langle v,v\rangle}$
we denote the Euclidean length of the vector $v$. If $G$ is an
additive group then $G_* = G\smin\{0\}$. By $\conv(C)$ we denote the
convex hull of a subset $C\ssq \R^d$. By $\textrm{Ex}(C)$ we denote
the extremal points of $\conv(C)$ and let $\conv_{\! \times}(C) =
\conv(C)\smin\textrm{Ex}(C)$.

If $v=\left(v_1,v_2\right) \in \R^2$ we let $v^\perp = \left(-v_2,v_1\right)$. Given $\lambda\neq 0$,
$v\in\R^2_*$ and $a\leq b \in \R\cup\{\pm \infty\}$ we let $\Llv =\lambda v +
\{v\}^\perp$ and
\begin{eqnarray*}
  C_v[a,b] & = & \{z\in\R^2 \mid a\langle z,v\rangle \leq
  \langle z,v^\perp \rangle \leq b\langle z,v \rangle \} \\
  \Llv[a,b] & = &  \{z\in \Llv \mid
  a\langle z,v\rangle \leq \langle z,v^\perp \rangle \leq b\langle z,v \rangle \}  \ .
\end{eqnarray*}
Note that thus $\Llv(a,b) = \Llv \cap C_v(a,b)$ and $C_v[a,a] = \R \cdot
(v+av^\perp)$. Further, we let
\begin{eqnarray*}
C^+_v[a,b] & = & \{ z\in C_v[a,b] \mid \langle z,v \rangle \geq 0\} \quad \textrm{and} \\
S_v[a,b] & = & \{ z\in\R^2 \mid \langle z,v\rangle \in [a,b]\} \ . 
\end{eqnarray*}
All these notions are used similarly for open and half-open intervals. 
\smallskip

The following basic observation is a direct consequence of the
definition in \eqref{e.rotation-subset}. Recall that
$\varphi_n(z)=(F^n(z)-z)/n$.
  \begin{lem} \label{l.convergence} Suppose $F$ is a lift of
   $f\in\homd$ and $U\ssq \T^d$. Then for all $\eps > 0$ there exists
    some $n_0=n_0(\eps) \in\N$ such that $\varphi_n(U) \ssq
    B_{\eps}(\rho_U(F)) \ \forall n \geq n_0$.
  \end{lem}
  The proof of the following statementis more or less identical to
  that of the connectedness of the rotation set in
  \cite{misiurewicz/ziemian:1989}, but we include it for the
  convenience of the reader.
\begin{lem} \label{l.connectedness} Let $F$ be a lift of
  $f\in\homeo_0(\T^d)$. For any $U\ssq \T^d$, the set $\rho_U(F)$ is
  compact. Further, if $U$ is connected then so is $\rho_U(F)$.
\end{lem} {\em Proof of Lemma~\ref{l.connectedness}.}  The fact that
$\rho_U(F)$ is compact follows immediately from the definition.
Suppose for a contradiction that $U$ is connected but $\rho_U(F)$ is
not. Then there exist disjoint open sets $V_1$ and $V_2$ with
$\rho_U(F) \ssq V_1 \cup V_2$ and $\rho_U(F) \cap V_i \neq \emptyset \
(i=1,2)$. Since $\rho_U(F)$ is compact, we may assume that $\eps =
d\left(V_1,V_2\right) > 0$. Lemma~\ref{l.convergence} implies that
there exists $n_0\in\N$ such that $\varphi_n(U) \ssq V_1\cup V_2 \
\forall n\geq n_0$. Further, since $\sup_{z\in\R^2} \|F(z)-z\| <
\infty$ there exists $n_1\in\N$ such that
\begin{equation}
  \left|\left(F^{n+1}(z)-z\right)/(n+1) - \left(F^n(z)-z\right)/n\right| \ <
  \ \eps \quad \forall n\geq n_1,\ z\in\R^2 \ . 
\end{equation}
It follows that if $n\geq n_1$ and $\varphi_n(U) \ssq V_i$, then $\varphi_k(U)
\ssq V_i \ \forall k\geq n$ and therefore $\rho_U(F) \ssq \overline{V_i}$. Since
this is not the case, $\varphi_n(U)$ must intersect both $V_1$ and $V_2$ for all
$n\geq n_1$. However, for $n\geq \max\{n_0,n_1\}$ we then obtain $\varphi_n(U)
\ssq V_1\cup V_2$ and $\varphi_n(U) \cap V_i\neq \emptyset \ (i=1,2)$. This
contradicts the connectedness of $\varphi_n(U)$. \qed

\medskip

\noindent
{\em Deviations from a constant rotation and invariant measures.}
For $f\in\homd$ with lift $F :\R^d \to \R^d$, $\rho\in\R^d$ and
$v\in\R^d_*$ we let 
\begin{equation}
  \label{e.deviations}
  D_n(z,\rho) \ :=  \ F^n(z)-z-n\rho \quad \textrm{and} \quad  D^v_n(z,\rho) \ := 
 \ \left\langle D_n(z,\rho)\right\rangle_v \ .  
\end{equation}
If we need to make the dependence on $f$ explicit, we also write
$D_{f,n}(z,\rho)$ and $D^v_{f,n}(z,\rho)$. 
For any $f$-invariant probability measure $\mu$, the rotation number of $f$ with
respect to $\mu$ is given by
\begin{equation}
  \label{e.mu-rotnum}
  \rho_\mu(F) \ = \ \int_{\T^d} F(z)-z \  d\mu(z) \ .
\end{equation}
When $F$ is fixed and no ambiguities can arise, we suppress it from the notation
and write $\rho_\mu$ instead of $\rho_\mu(F)$. By $\supp(\mu)$ we denote the
topological support of $\mu$.
\begin{lem}
  \label{l.measure-drift} Suppose $F : R^d \to \R^d$ is a lift of $f\in\homd$
  and $\mu$ is an ergodic $f$-invariant probability measure.  Then
 there exists no constant $s>0$ with the property that for $\mu$-a.e.\
  $z\in\R^2$ there is a positive integer $n_z$ such that $D^v_{n_z}(z,\rho_\mu) \geq
  sn_z$.
\end{lem}
\myproof 
We suppose for a contradiction that a constant $s>0$ with the above
property exists. We fix an $f$-invariant set $\Omega \ssq \torus$ of
measure $\mu(\Omega) = 1$ such that for all $z\in\Omega$ there exists
$n_z\in\N$ with
\begin{equation} \label{e.dev-nz} D^v_{n_z}(z,\rho_\mu) \ \geq \ sn_z \ .
\end{equation}
In addition, we assume that 
\begin{equation} \label{e.rho-mu-convergence} \nLim \left(F^n(z)-z\right)/n \ = \ \rho_\mu
  \quad \forall z\in\Omega \ .
\end{equation}
Given any $z_0 \in \Omega$, we recursively define a sequence of integers $n_i$
by $n_0 = 0$ and $n_{i+1} = n_i + n_{F^{n_i}(z_0)}$. Then we obtain
\begin{eqnarray*}
  D^v_{n_k}(z_0,\rho_\mu) & = & \left\langle F^{n_k}(z_0) - z_0 - n_k\rho_\mu \right\rangle_v \\
  & = &
  \left\langle \sum_{i=0}^{k-1} F^{n_{i+1}}(z_0) - F^{n_{i}}(z_0) - (n_{i+1}-n_{i})\rho_\mu\right\rangle_v 
  \\ & = &  \left\langle \sum_{i=0}^{k-1} F^{n_{F^{n_{i}}(z_0)}}\left(F^{n_i}(z_0)\right) - F^{n_i}(z_0) -
    (n_{F^{n_i}(z_0)})\rho_\mu\right\rangle_v  \\ & \stackrel{(\ref{e.dev-nz})}{\geq} &
  \sum_{i=0}^{k-1} s(n_{i+1}-n_i) \ = \ sn_k  \ . 
\end{eqnarray*}
Hence $\kLim D^v_{n_k}(z_0,\rho_\mu)/n_k \geq s$, contradicting
(\ref{e.rho-mu-convergence}) which implies $\nLim D^v_n(z_0,\rho_\mu) = 0$.
\qed\medskip
\pagebreak

\noindent
{\em A reduction lemma.}
Any integer matrix $M \in \textrm{GL}(d\times d,\Z)$ induces a toral
endomorphism $g_M : \T^d \to \T^d,\ \pi(z) \mapsto \pi(Mz)$, and $g_M$ is
invertible if and only if $M \in \textrm{SL}(d,\Z)$. The following lemma
describes how a coordinate transformation by such a map $g_M$ acts on the
rotation set. 

\begin{lem}[Reduction Lemma] \label{l.reduction} Suppose $f\in\homeo_0(\T^d)$
  has lift $F:\R^d \to \R^d$, $U\ssq \T^d$ and $M \in \mathrm{SL}(d,\Z)$. Then
  the following hold.  \alphlist
\item Let $\tilde f = g_M^{-1} \circ f \circ g_M \in
  \homeo_0(\T^d)$ with lift $\tilde F = M^{-1} \circ F \circ M$. Then
  \begin{equation}
    \label{e.rotset-transformation}
    \rho_{g_M^{-1}U}(\tilde F) \ = \ M^{-1}(\rho_U(F)) \ . 
  \end{equation}
\item If $\rho_U(F) \ssq \Llv$ then $\rho_{g_M^{-1}U}(\tilde F) \ssq
  M^{-1}(\Llv) = L_{\tilde \lambda,\tilde v}$, where $\tilde v = M^tv$ and
  $\tilde \lambda = \lambda \|v\|^2/\|\tilde v\|^2$. Further $D^v_{f,n}(z,\rho)
  = D^{\tilde v}_{\tilde f,n}(M^{-1}z,M^{-1}\rho) \ \forall n\in\N,\
  z,\rho\in\R^d$.
\item Let $1 \leq k < d$ and suppose that $w_1 \ld w_k$ are linearly independent
  integer vectors and $\conv_{\! \times}(\{w_1\ld w_k\})$ contains no further
  integer vectors. (If $k=1$ this just means that the entries of $w_1$ are
  relatively prime.) Then there exist integer vectors $w_{k+1} \ld w_d$ such
  that $\det(w_1\ld w_d) = 1$. \listend
\end{lem}\noindent
Note that for any integer vector $w \in \Z^d$ with
relatively prime entries part (c) allows to perform a linear
coordinate transformation on $\T^d$ such that $w$ becomes a base
vector.\medskip

\myproof \alphlist
\item Suppose $z_i \in \pi^{-1}U,\ n_i \nearrow \infty$ and $\iLim
  \left(F^{n_i}\left(z_i\right)-z_i\right)/n_i = \rho$. Then
  \begin{eqnarray*}
    M^{-1}\rho & = & M^{-1} \left(\iLim \left(F^{n_i}\left(z_i\right)-z_i\right)/n_i\right)  
    \\ & = & \iLim \left(M^{-1}\circ F^{n_i} \circ M \left(M^{-1}z_i\right)-M^{-1}z_i\right)/n_i \\
    & = & \iLim \left(\tilde F^n\left(M^{-1}z_i\right)-M^{-1}z_i\right)/n_i \ 
    \in \ \rho_{M^{-1}U}(\tilde F) \ .
  \end{eqnarray*}
  This shows that $M^{-1}\left(\rho_U\left(F\right)\right) \subseteq
  \rho_{g_M^{-1}U}\left(\tilde F\right)$ and since $f = g_M \circ \tilde f \circ
  g_M^{-1}$ the opposite inclusion follows in the same way.
\item $\rho_{g_M^{-1}U}(\tilde F) \ssq M^{-1}(\Llv)$ holds by part (a). Further,
  we have $M^{-1}(\Llv) = M^{-1}(\lambda v) + M^{-1}\{v\}^\perp$. Since $\tilde
  v = M^t v \perp M^{-1}\{v\}^\perp$ if follows that $M^{-1}(\Llv) =
  L_{\tilde\lambda,\tilde v}$ for some $\tilde\lambda \in \R$ and we have
  $\tilde\lambda \|\tilde v\| = \langle M^{-1}(\lambda v),\tilde v/\|\tilde
  v\|\rangle = (\lambda/\|\tilde v\|) \cdot \langle v,v \rangle$. Finally, in
  order to check that $D_{f,n}^v(z,\rho) = D_{\tilde f,n}^{\tilde
    v}(M^{-1}z,M^{-1}\rho)$ let $z\in \R^d$. Then
  \begin{eqnarray*}
    D_{\tilde f,n}^{\tilde v}(M^{-1}z,M^{-1}\rho) & = & \langle \tilde F^n(M^{-1}z) - 
    M^{-1}z - nM^{-1}\rho,\tilde v\rangle  \\
    & = & \langle M^{-1} \circ F^n (z) - M^{-1}z-nM^{-1}\rho,M^tv\rangle  \\
    & = & \langle F^n(z)-z-n\rho,v\rangle  \ = \ D_{f,n}^v(z,\rho) \ . 
  \end{eqnarray*}
\item Choose integer vectors $w_{k+1} \ld w_d$ such that $\spann(w_1 \ld w_d) =
  \R^d$. Then $|\det(w_1 \ld w_d)| \geq 2$ if and only if $\conv_{\! \times}(w_1
  \ld w_d)$ contains an integer vector. In this case we replace one of the
  vectors $w_{k+1} \ld w_d$ by an integer vector in $\conv_{\! \times}(w_1 \ld
  w_d)$ such that the new set of vectors still spans $\R^d$. This reduces the
  absolute value of the determinant, and after a finite number of steps we
  arrive at $\det(w_1\ld w_d) = \pm 1$. Replacing $w_d$ by $-w_d$ if necessary
  we obtain $\det(w_1\ld w_d) = 1$.  \qed \listend

\section{Rotation subsets on the universal cover}
\label{Semilocal-rotsets} Throughout this section, we suppose that $G$
is the lift of a toral homeomorphism $g\in\homtwo$, $\widehat{U} \ssq
\R^2$ is bounded and connected and $G(\widehat{U}) \cap \widehat{U}
\neq \emptyset$. Further, we assume that $\lambda\neq 0$ and
$v\in\R^2_*$. 

In order to control the whole rotation set by using assumptions on
$\rho_{\widehat U}(G)$, we proceed in several steps. The first is to obtain some
information about the extremal points of the rotation set.
\begin{lem} \label{l.extremalpoints} Suppose $\rho_{\widehat{U}}(G) =
  L_{\lambda,v}[a,b]$. Then all extremal points of $\rho(G)$ belong to
  $C_v[a,b]$.
\end{lem}
\myproof Suppose that $\lambda > 0$. (Otherwise we replace $v$ by $-v$.)
Performing a linear change of coordinates via Lemma~\ref{l.reduction} if
necessary we may assume that
\begin{equation} \label{e.Llv-proj}
\pi_1(\Llv[a,b]) \ \ssq \ (0,\infty) \ ,
\end{equation}
such that in particular $C^+_v[a,b]\smin \{0\} \ssq (0,\infty) \times
\R$. Further, we may assume that $\widehat{U}$ intersects $\{0\}\times \R$,
otherwise we replace it by an integer translate and/or one of its iterates.%
\foot{Note that by assumption $\bigcup_{n\in\N_0} G^n(\widehat U)$ is connected
  and $\pi_1\circ G^n(\widehat{U})$ goes to $\infty$ as $n\to\infty$ due to
  (\ref{e.Llv-proj}). Hence, one of the iterates of $\widehat U$ has to
  intersect an integer vertical $\{m\}\times\R$.}
Let $V := \bigcup_{n\in\N_0} G^n(\widehat{U})$.  As $G(\widehat{U})
\cap \widehat{U} \neq \emptyset$, the set $V$ is connected. We claim
that for sufficiently large $l\in\N$ the integer translate $V-(0,l)$
is disjoint from $V$. In order to see this let $r:=\sup_{z\in
  \widehat{U}} \|z\|$. Due to (\ref{e.Llv-proj}) and
Lemma~\ref{l.convergence} only a finite number of iterates of
$\widehat{U}$ intersect $[-r,r]\times \R$. Therefore $V$ and
$\widehat{U}-(0,l)$ are disjoint for large $l$.  Hence, if the orbit
of $\widehat{U}-(0,l)$ intersects $V$ then it must first intersect
$\widehat{U}$. However, by the same argument the orbit of
$\widehat{U}-(0,l)$ can only intersect a finite number of its vertical
integer translates, such that for sufficiently large $l$ we have
$V\cap [V-(0,l)] = \emptyset$ as required.

Now let $y_1=\inf\{y\in\R \mid (0,y) \in V-(0,l)\}$, $y_2=\sup\{y\in\R \mid
(0,y) \in V\}$ and define $W$ as the union of $V$, $V-(0,l)$ and the vertical
arc from $(0,y_1)$ to $(0,y_2)$. Let $Y$ be the unique connected component of
$\R^2\smin \overline{W}$ which is unbounded to the left and $A = \R^2 \smin
Y$. 

\begin{figure}[h]
  \begin{center} \epsfig{file=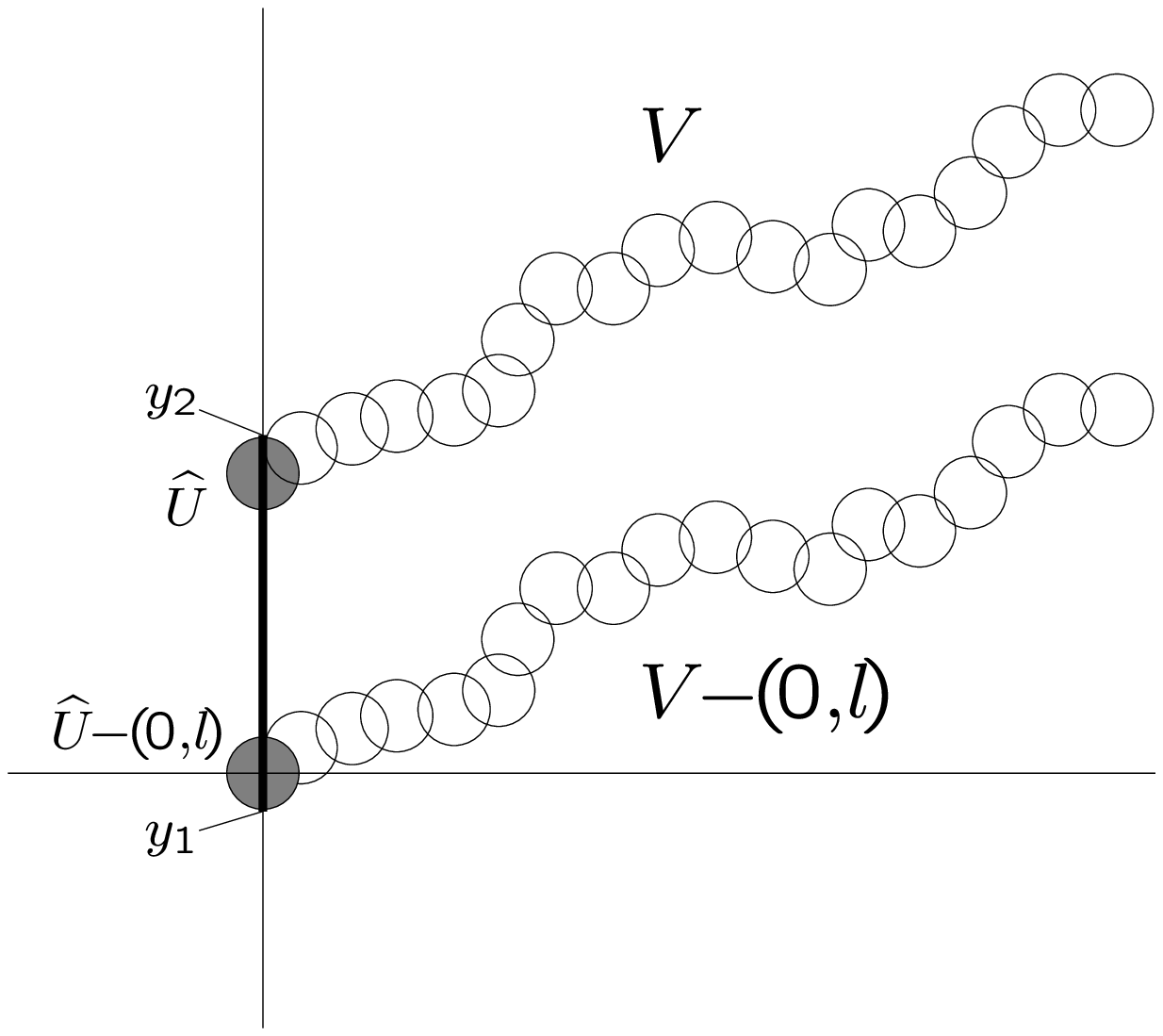, clip=,
      width=0.4018\linewidth} \hspace{1eM} \epsfig{file=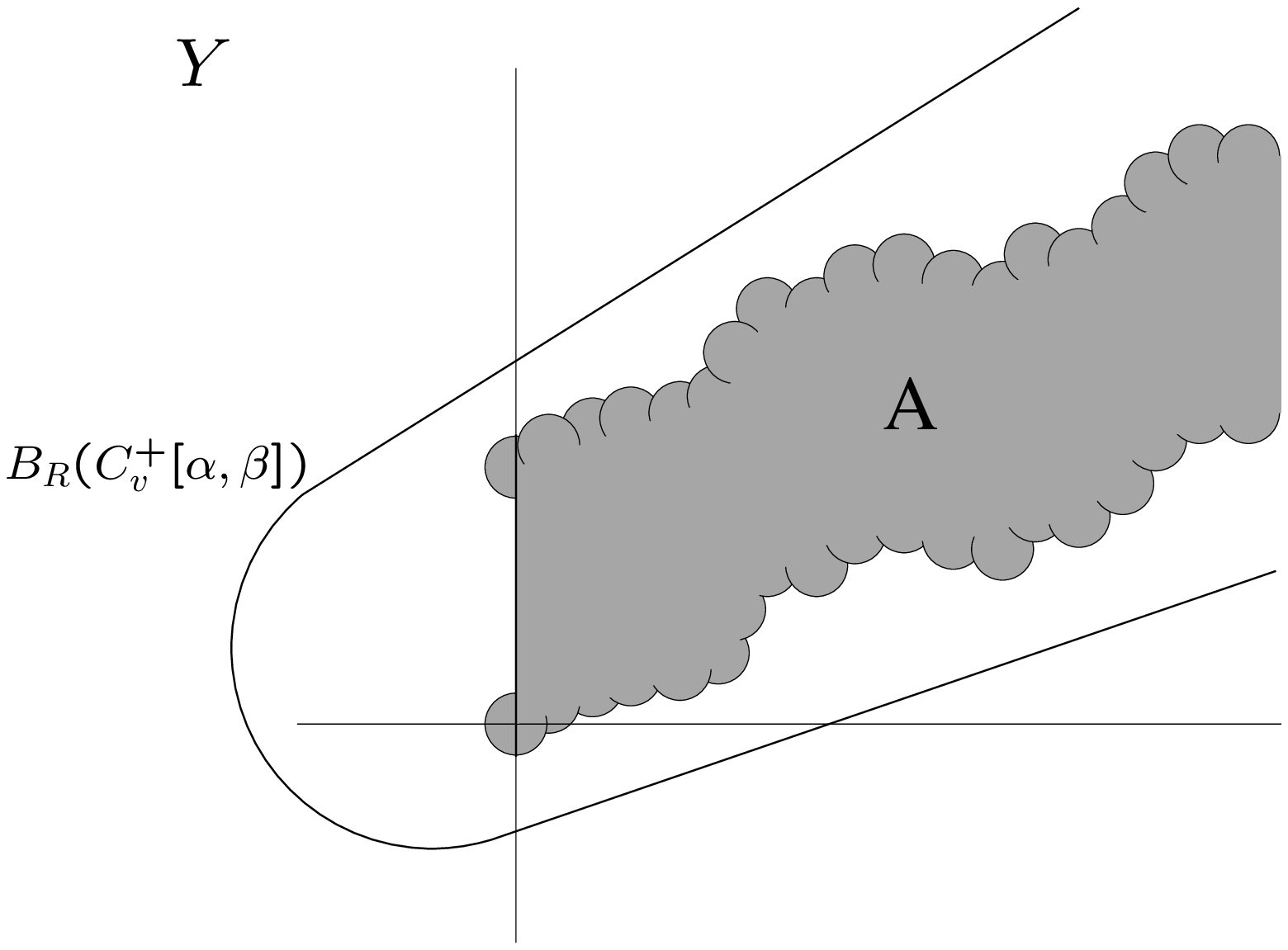, clip=,
      width=0.539\linewidth}\end{center}
  \caption{\small Construction of the sets $W$ on the left and $A$ on the
    right. \label{f.2} }
\end{figure}

The following three remarks about these objects will be helpful.
First, as $\rho_{\widehat U}(G) \ssq \Llv[a,b]$ and $\lambda>0$, we
have that for all $\alpha<a$ and $\beta>b$ there exists a constant
$R=R(\alpha,\beta)>0$ such that
\begin{equation} \label{e.A-inclusion}
  A \ \ssq \ B_R(C^+_v(\alpha,\beta)) \ . 
\end{equation}
Secondly, due to the definition of $W$, its connectedness and
(\ref{e.Llv-proj}), the set $Z := (\R^+ \times \R) \smin A = (\R^+ \times \R)
\cap Y$ consists of exactly two connected components. These can be defined as
follows. Fix any $\zeta_0 \in Y$ with $\pi_1(\zeta_0) < 0$. For any $\zeta \in Z$,
there is a path $\gamma_\zeta$ in $Y$ from $\zeta$ to $\zeta_0$. Let $y_\zeta$ be
the second coordinate of the first point in which $\gamma_\zeta$ intersects the
vertical axis. The fact whether $y_\zeta$ lies below $y_1$ or above $y_2$ does
not depend on the choice of the path, since this would contradict the
connectedness of $W$. Hence, $Z^- = \{\zeta \in Z \mid y_\zeta < y_1\}$ and $Z^+
= \{ \zeta \in Z \mid y_\zeta > y_2\}$ form a partition of $Z$ into two
connected components.

Thirdly, there holds $\R^+ \times \R \ssq \kcup A + (0,kl)$. Consequently, for
any $m\in \N$ the set $A \cap \pi_1^{-1}[m,m+1)$ contains a fundamental domain
of $\torus$, that is, $\pi(A \cap \pi_1^{-1}[m,m+1)) = \torus$.\medskip

It is important to note, however, that $A$ is not $G$-invariant. Yet, in order
to obtain control over the full rotation set via $A$ we will need to ensure that
orbits `moving to the right' become `trapped' in $A$ (or one of its integer
translates). Hence, the following statement is crucial for our purposes.
  \begin{claim} \label{c.rho} There exists a constant $K>0$ such that $z\in A
    \cap \pi_1^{-1}[K,\infty)$ implies $G^{\pm 1}(z) \in A$.
\end{claim}
\myproof We show that there exists a constant $K' > 0$ such that $z \in Z \cap
\pi_1^{-1}[K',\infty) $ implies $G^{\pm 1}(z) \in Z$. If we let
\begin{equation} \label{e.M}
M\ :=\ \sup_{z\in\R^2}\|G(z)-z\| \ = \ \sup_{n\in\R^2} \|G^{-1}(z)-z\| \ ,
\end{equation}
then for any $z\in\pi_1^{-1}[K'+M,\infty)$ this means that $G^{\pm 1}(z)
\in Z$ implies $z\in Z$ and hence $z\in A$ implies $G^{\pm 1}(z) \in A$. Thus we
can choose $K=K'+M$.

From (\ref{e.Llv-proj}) and Lemma~\ref{l.convergence} we deduce that there
exists $n_0\in\N$ such that $\pi_1\circ G^n(\widehat{U}) \ssq (4M,\infty) \
\forall n\geq n_0$. Let $K'>0$ such that for all $j\leq n_0$ there holds
$\pi_1\circ G^j(\widehat{U}) \ssq [0,K')$. Then we have
\begin{equation}
  \label{e.noreturn}
  \pi_1\circ G^n(\widehat{U}) \cap [K',\infty) \neq \emptyset \quad
  \follows \quad \pi_1\circ G^k(\widehat{U}) \cap [0,4M] = \emptyset \quad  \forall k\geq n \ .
\end{equation}
The same statement applies to $\widehat{U}+(0,l)$.  Due to (\ref{e.A-inclusion}) there
exists $C>0$ such that
\[
B \ := \ [0,4M] \times [C-M,\infty) \ \ssq \ Z^+ \ .
\]
Let $z^*=(3M,C)$ and fix $z \in Z^+$ with $\pi_1(z) \geq K'$. Then, since
$Z^+$ is open and connected, there is a simple path $\gamma : [0,1] \to Z^+$
from $z$ to $z^*$.  We claim that $\gamma$ can be chosen such that its image
is contained in $Z^+ \cap \pi^{-1}[3M,\infty)$.  Suppose not and let $t_0 :=
\min\{t\in[0,1] \mid \gamma(t) \in B\}$ and $\Gamma = \{\gamma(t) \mid
t\in[0,t_0]\}$. Then $\Gamma$ divides the set $([0,\pi_1(z)] \times \R) \smin B$
into exactly two connected components $D^+$ and $D^-$ that are unbounded above,
respectively below.
\begin{figure}[h]
  \begin{center} \epsfig{file=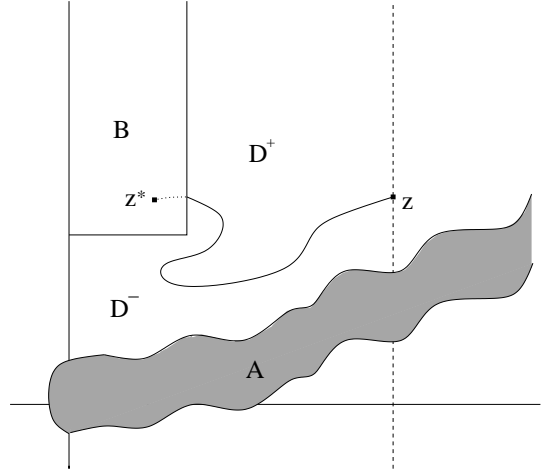, clip=,
      width=0.5\linewidth} \end{center}
  \caption{\small The domains $D^-$ and $D^+$. \label{f.4} }
\end{figure}
Now, if $\overline{W}$ does not intersect $D^+ \cap \pi^{-1}_1[0,4M)$ then
$D^+\cap \pi^{-1}_1[0,4M) \ssq Z^+$, and it is easy to see that in this case
either $\gamma$ does not intersect $\pi_1^{-1}[0,3M)$ or we can modify it to
that end. Otherwise, there must be some $k\in\N$ such that $G^k(\widehat{U})$ or
$G^k(\widehat{U})-(0,l)$ intersects $D^+\cap \pi^{-1}_1[0,4M)$. However, since
$\widehat{U}$ intersects $D^-$ and the set $\bigcup_{i=0}^k G^i(\widehat{U})$ is
connected, this implies that there must be some $n\leq k$ such that $\pi_1\circ
G^n(\widehat{U})$ intersects $[\pi_1(z),\infty) \ssq [K',\infty)$. This
contradicts (\ref{e.noreturn}).

Summarising, we have found a path $\gamma$ from $z$ to $z^*$ which is contained
in $Z^+ \cap \pi^{-1}_1[3M,\infty)$. In particular, $\gamma$ is contained in the
complement of $\overline W$. Consequently, the path $G\circ \gamma$ is contained
in the complement of $G(\overline{W})$. At the same time, it is also contained
in $\pi^{-1}_1[2M,\infty)$. However, it follows from the construction of $W$ and
the definition of $M$ that $G(\overline W) \cap \pi^{-1}_1[2M,\infty) =
\overline W \cap \pi^{-1}_1[2M,\infty)$. Hence, the path $G \circ \gamma$ is
contained in the complement of $\overline{W}$ as well.  Furthermore, it joins
$G(z)$ to the point $G(z^*)$. Since the latter is contained in $B \ssq Z^+$,
this implies that $ G(z)$ is equally contained in $Z^+$. When $z\in Z^-$ the
argument is similar. In the same way one can show that $G^{-1}(z) \in Z$, and
this proves the claim.  \roundqed \medskip

In order to complete the proof of Lemma~\ref{l.extremalpoints},
suppose that $\rho$ is an extremal point of $\rho(G)$ which is not
contained in $C_v[a,b]$. By performing a linear change of coordinates
again if necessary, we may assume that $\pi_1(\rho) > 0$.\foot{Choose
  a basis of integer vectors $w_1,w_2$ with $\det(w_1,w_2)=1$ such
  that both $w_1$ and $\rho$ lie to the right of the oriented line $\R
  w_2$ and $w_2 \notin C_v[a,b]$ (the latter ensures that
  (\ref{e.Llv-proj}) remains valid.) Then apply
  Lemma~\ref{l.reduction}. } Since $\rho$ is realised by an ergodic
invariant measure \cite[Corollary~3.5]{misiurewicz/ziemian:1989},
there exists a point $z_0\in\R^2$ with
\begin{equation} \label{e.z_0-convergence}
\nLim (G^n(z_0)-z_0)/n \ = \ \rho\ . 
\end{equation}
Let $z_n = G^n(z_0)$ and ${\cal O}^+(z_0) = \{z_n \mid n\geq 0\}$.
Then, due to (\ref{e.z_0-convergence}), for every $\gamma< 0 < \delta$
there exists $\tilde R>0$ such that ${\cal O}^+(z_0) \ssq z_0 +
B_{\tilde R}(C^+_\rho[\gamma,\delta])$. If $\gamma,\delta$ are chosen
sufficiently close to 0 and $\alpha,\beta$ in (\ref{e.A-inclusion})
are sufficiently close to $a$ and $b$ then $B_R(C_v[\alpha,\beta])
\cap \left(z_0+B_{\tilde R}(C_\rho[\gamma,\delta])\right)$ is bounded
(recall that $\rho$ is not contained in $C_v[a,b]$). Further, by
replacing $z_0$ with an integer translate if necessary, we may assume
that $z_0 \in A$ and $\pi_1(z_0) \geq K + \tilde R$, where $K$ is
chosen as in Claim~\ref{c.rho}.  If follows that $\pi_1(z_n) \geq K$
and hence $z_n \in A \ \forall n\in\N$. Consequently ${\cal O}^+(z_0)
\ssq B_R(C_v[\alpha,\beta]) \cap \left(z_0+B_{\tilde
    R}(C_\rho[\gamma,\delta])\right)$ such that ${\cal O}^+(z_0)$ is
bounded, contradicting (\ref{e.z_0-convergence}).  Hence, all extremal
points of $\rho(G)$ must be contained in $C_v[a,b]$.  \qed \medskip


In the opposite way, information about the extremal points of $\rho(G)$ allows
to draw conclusions about the behaviour of the iterates of $\widehat U$. 
\begin{lem}
  \label{l.cone-intersection} Suppose that $\rho_{\widehat{U}}(G) =
  L_{\lambda,v}[a,b]$, $\gamma \in [a,b]$ and $\rho \in C_v[\gamma,\gamma]$ is
  an extremal point of $\rho(G)$. Then given any $\eps > 0$ there exists
  $N\in\N$ such that $G^n(\widehat{U}) \cap C^+_v(\gamma-\eps,\gamma+\eps) \neq
  \emptyset \ \forall n\geq N$.
\end{lem}
\myproof As in the proof of Lemma~\ref{l.extremalpoints} we assume $\lambda > 0,
\pi_1(v) > 0$ and $C^+_v[a,b]\smin\{0\} \ssq (0,\infty) \times \R$. Let $l$ be
the integer in the definition of the set $A$ above and define $\widehat{U}' :=
\widehat{U} \cup [\widehat{U}-(0,l)]$. Suppose for a contradiction that
$G^n(\widehat{U}) \cap C^+_v(\gamma-\eps,\gamma+\eps) = \emptyset$ for
infinitely many $n\in \N$.  Slightly reducing $\eps$ if necessary, we may assume
that $G^n(\widehat{U}') \cap C^+_v(\gamma-\eps,\gamma+\eps) = \emptyset$ for
infinitely many $n\in \N$.
\begin{figure}[h]
  \begin{center} \epsfig{file=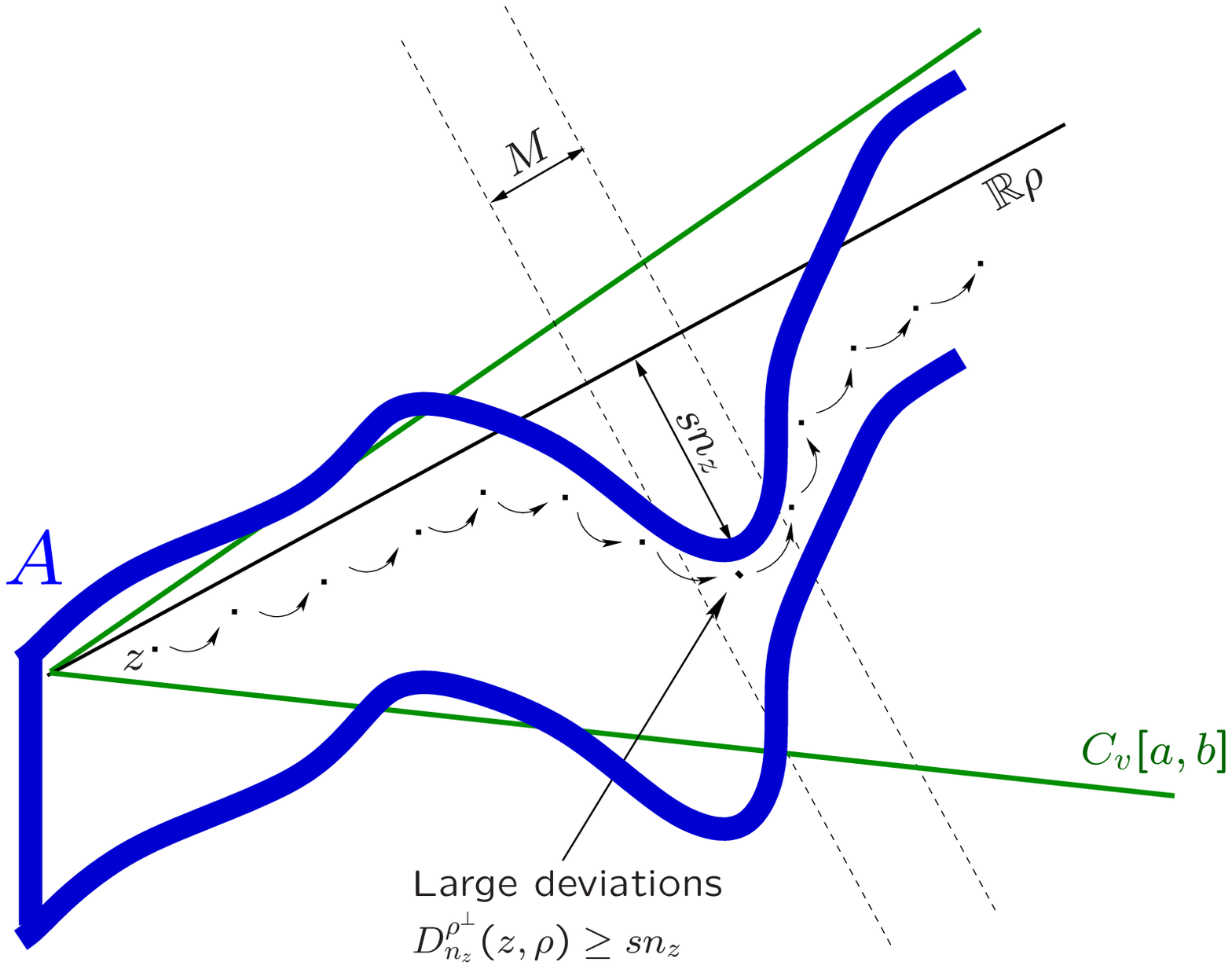, clip=,
      width=0.6\linewidth} \end{center}
  \caption{\small Strategy for the proof of Lemma~\ref{l.cone-intersection}: The
    existence of large 'lacunae' in the set $A$ (with blue boundary) forces
    orbits inside $A$ to deviate far from the line $\R\rho$. In particular, this
    is true for almost all orbits w.r.t.\ the measure $\mu$ realising the
    rotation vector $\rho \in \textrm{Ex}(\rho(G))$. This leads to a
    contradiction with Lemma~\ref{l.measure-drift}. \label{f.9} }
\end{figure}
We first consider the case where $\vbr{\rho} > 0$, such that $\rho \in
C^+_v[a,b]$ by Lemma~\ref{l.extremalpoints}. 

Due to Lemma~\ref{l.convergence} the fact that $\rho_{\widehat{U}}(G)
\ssq \Llv$ implies that $\vbr{G^nz}/n$ converges uniformly to
$\lambda$ on $\widehat{U}'$ as $n\to\infty$.  Hence, for any
$\delta>0$ there exists $N(\delta)\in\N$ such that
\begin{equation}
  \label{e.N_delta} G^n(\widehat{U}')  \ \ssq \ S_v[(1-\delta)n\lambda,(1+\delta)n\lambda] 
  \quad \forall n\geq N(\delta) \ . 
\end{equation}
As $\pi_1(v) > 0$, this implies that $\inf \pi_1(G^n(\widehat{U}')) \to \infty$
as $n\to\infty$. Therefore (\ref{e.A-inclusion}) yields that for given $\alpha<
a$ and $\beta>b$ and sufficiently large $n$ there holds $G^n(\widehat{U}') \ssq
C^+_v[\alpha,\beta]$ (apply (\ref{e.A-inclusion}) with $\tilde \alpha \in
(\alpha,a)$ and $\tilde \beta \in (b,\beta)$ to get rid of the constant
$R$). Consequently, if $n$ is large then $G^n(\widehat{U}') \cap
C^+_v(\gamma-\eps,\gamma+\eps)=\emptyset$ implies $G^n(\widehat{U}')\ssq
C^+_v[\alpha,\gamma-\eps]$ or $G^n(\widehat{U}') \ssq
C^+_v[\gamma+\eps,\beta]$. We assume that $G^n(\widehat{U}') \ssq
C^+_v[\alpha,\gamma-\eps]$ for infinitely many $n \in\N$, the other case is
symmetric.

Since $\rho \in \textrm{Ex}(\rho(G))$, there exists an ergodic measure $\mu$
with $\rho_\mu(F) = \rho$ \cite[Corollary~3.5]{misiurewicz/ziemian:1989}. Let
$\Omega \ssq \torus$ be such that $g(\Omega)=\Omega$, $\mu(\Omega)=1$ and
\begin{equation}
  \label{e.mu-convergence} \nLim (G^n(z)-z)/n \ = \ \rho \quad \forall z\in\pi^{-1}(\Omega) \ . 
\end{equation}
We will show that, in contradiction to Lemma~\ref{l.measure-drift},
for some $s>0$ there holds
\begin{equation}\label{e.s-contradiction}
\forall z\in\Omega \ \exists n_z\in\N : \quad D^{-\rho^\perp}_{n_z}(z,\rho) \ \geq \  s n_z \ .
\end{equation}
In order to do so, fix $z\in\Omega$. Let $A$ be as in the proof of
Lemma~\ref{l.extremalpoints}. Due to (\ref{e.mu-convergence}) there
exists a lift $z_0\in A$ of $z$ such that $\pi_1(z_n) \geq K \ \forall
n\in\N$, where $z_n = G^n(z_0)$ and $K$ is chosen as in Claim~\ref{c.rho}.
Consequently Claim~\ref{c.rho} implies that
\begin{equation}
  \label{e.z_n-in-A}
  z_n \ \in \ A \quad \forall n\in\N \ . 
\end{equation}
Due to (\ref{e.N_delta}) and the fact that $M$
defined in (\ref{e.M}) is finite, there exists $\eta > 0$ and $N_1\in\N$ 
such that for any $n\geq N_1$ there holds
\begin{equation}
  \label{e.cone-inclusion} G^n(\widehat{U}')  \ssq  C_v^+[\alpha,\gamma-\eps] \quad \follows 
  \quad G^k(\widehat{U}')  \ssq  C_v^+[\alpha,\gamma-\eps/2] \qquad  \forall k\in[(1-\eta)n,(1+\eta)n]  \ .
\end{equation}
Now, choose $\delta$ in (\ref{e.N_delta}) sufficiently small, such that
$(1+\delta)(1-\eta) < 1 < (1-\delta)(1+\eta)$. Then choose $N_2\in\N$ such that
for all $n\geq N_2$ there holds
\begin{eqnarray} \label{e.interval-inclusion} [\lambda n-M,\lambda n+M] & \ssq &
  [(1+\delta)(1-\eta)\lambda n,(1-\delta)(1+\eta)\lambda n] \ ,\\
  \lambda n-M & \geq & M\cdot N(\delta)\cdot \|v\|+\sup \left\langle\widehat U'\right\rangle_v \
  . \label{e.geq-N_delta}
\end{eqnarray}
Suppose $n\geq N_3:=\max\{N(\delta)/(1-\delta),N_1,N_2\}$ and $G^n(\widehat{U}) \ssq
C_v^+[\alpha,\gamma-\eps]$.  Then by combining (\ref{e.N_delta}) and
(\ref{e.cone-inclusion})--(\ref{e.geq-N_delta}) we obtain that 
\begin{equation}
  G^k(\widehat{U}')
  \cap S_v(\lambda n-M,\lambda n+M) \ \ssq \ C^+_v[\alpha,\gamma-\eps/2] \quad
  \forall k \in\N \ .
\end{equation}
(Treat the cases $k\leq N(\delta),\ k\in(N(\delta),(1-\eta)n)$ and $k\geq
(1+\eta)n$ separately to show that for all such $k$ the set $G^k(\widehat U)$
does not intersect $S_v(\lambda n-M,\lambda n+M)$ and then use
(\ref{e.cone-inclusion}) for the remaining~$k$.) This means in particular that
\begin{equation} \label{e.A-pocket} A\cap S_v(\lambda n-M,\lambda n+M) \ \ssq \
  C_v^+[\alpha,\gamma-\eps/2] \quad \quad \forall n\geq N_3 : G^n(\widehat{U}) \ssq
  C^+_v[\alpha,\gamma-\eps] \ .
\end{equation}
Now, as $\R\rho=C_v[\gamma,\gamma]$ exist constants $r>0$ and $N_4\geq
N_3$ such that
\begin{equation} \label{e.r-deviationbound}
  \langle z'-z_0\rangle_{\rho^\perp} \leq -rn \quad \forall n\geq N_4,\ 
  z'\in S_v(\lambda n-M,\lambda n+M) \cap C^+_v[\alpha,\gamma-\eps/2] \ .
\end{equation}
Further, since $z$ and its lift $z_0$ are fixed and due to
(\ref{e.mu-convergence}) (applied to $z_0$) there exist constants $N_5=N_5(z)
\geq N_4$ and $c>0$, with $c$ only depending on $\rho,\gamma$ and $\eps$, such
that
\begin{equation} \label{e.nz}
  \forall n\geq N_5 \ \exists n_z \leq cn : \quad z_{n_z} \in S_v(\lambda n-M,\lambda
  n+M)
\end{equation} 
(Note that the orbit of $z_0$ has to pass through $S_v[\lambda
n-M,\lambda n+M]$ by definition of $M$, and due to
(\ref{e.mu-convergence}) this happens approximately at time
$n\lambda/\vbr{\rho}$.)  Choose $n\geq N_5$ with $G^n(\widehat U)$
$\ssq$ $C_v^+[\alpha,\gamma-\eps]$. Let $n_z$ be as in (\ref{e.nz}).
Since $z_{n_z} \in A$ by (\ref{e.z_n-in-A}), we obtain $\langle
z_{n_z}-z_0\rangle_{\rho^\perp} \leq -rn$ from (\ref{e.A-pocket}) and
(\ref{e.r-deviationbound}). Thus, if $s=r/c$ then
$D_{n_z}^{-\rho^\perp}(z,\rho) \geq sn_z$. Since $z\in\Omega$ was
arbitrary, this proves (\ref{e.s-contradiction}).

Finally, if $\vbr{\rho} < 0$ then we can proceed in the same way by regarding
the inverse of $G$. Due to the symmetry in the statement of Claim~\ref{c.rho} we
obtain that for a suitable lift $z_0$ of $z\in\Omega$ the whole backwards orbit
of $z_0$ remains in $A$, and the remaining argument is exactly the same as
before.  \qed \medskip

\begin{cor}
  \label{c.linear-spreading} Suppose $\rho_{\widehat{U}}(G) = L_{\lambda,v}[a,b]$ with
  $a<b$. Then there exist constants $c>0$ and $N'\in\N$ such that
  \begin{equation}
    \label{e.linear-spreading} \sup_{z\in \widehat{U}} \vpbr{G^n(z)} - \inf_{z\in \widehat{U}}\vpbr{G^n(z)} 
    \ > \ cn \quad \forall n \geq N' \ . 
  \end{equation}
\end{cor}
\myproof As $L_{\lambda,v}[a,b] = \rho_{\widehat{U}}(G) \ssq \rho(G)$,
the set $\rho(G)$ must have at least two linearly independent extremal
points $\rho_1,\rho_2\neq 0$. Due to Lemma~\ref{l.extremalpoints}
these are contained in $C_v[a,b]$, such that $\rho_1\ssq
C_v[\gamma_1,\gamma_1]$ and $\rho_2\in C_v[\gamma_2,\gamma_2]$ for
some $\gamma_1\neq\gamma_2\in\R$. The statement now follows from
Lemma~\ref{l.cone-intersection} together with the fact that
$\inf\langle G^n(\widehat U)\rangle_v \to \infty$ as $n\to\infty$.  \qed \medskip

We can now
describe two situations in which the rotation subset of $\widehat U$ determines
the whole rotation set completely, or at least forces it to be contained in a
line segment.
\begin{lem}
  \label{l.uniform-v-speed} Suppose $\rho_{\widehat{U}}(G) = L_{\lambda,v}[a,b]$
  with $a<b$. Then $\rho(G) = L_{\lambda,v}[a,b]$.
\end{lem}
\myproof Due to Lemma~\ref{l.extremalpoints}, it suffices to show $\rho(G) \ssq
\Llv$. Suppose for a contradiction that $\rho(G) \nsubseteq \Llv$. Then there
exists an extremal point $\rho \in \textrm{Ex}(\rho(G)) \smin \Llv$. We assume
w.l.o.g.\ that $\|v\| = 1$ and $\vbr{\rho} > \lambda$. Since $\rho$ is realised
by an ergodic measure, there exists $z_0 \in \R^2$ with $\nLim (z_n-z_0)/n =
\rho$, where $z_n = G^n(z_0)$ as above. Fix $\eta > 0$ and $k_0\in\N$ such that
\begin{equation}
  \label{e.k_0-unifspeed} \vbr{z_k-z_0} \ > \ (1+\eta)\lambda k \quad \forall k \geq k_0 \ . 
\end{equation}
Further, fix $\delta > 0$ such that 
\begin{equation}
  \label{e.delta-unifspeed}
  \delta\left(1+15M/c\right) \ < \ \eta \ , 
\end{equation}
with $M$ defined by (\ref{e.M}) and $c$ as in
Corollary~\ref{c.linear-spreading}. Choose $n_0\in\N$ such that
\begin{eqnarray}
  \label{e.n_0-unifspeed}
  G^n(\widehat{U}) & \ssq & S_v[(1-\delta)\lambda n,(1+\delta)\lambda n] \quad 
 \forall n\geq n_0 \ \quad \textrm{and} \\
  \delta \lambda n_0 & \geq & 2 \ . \label{e.n_0-second}
\end{eqnarray}
Then choose $k\geq k_0$ and $n\geq n_0$ such that
\begin{equation}
  \label{e.k_0-n_0-ineq}
  (4M+2)k \ \leq \ cn \ \leq 5Mk \ .
\end{equation}

\begin{figure}[h]
  \begin{center} \epsfig{file=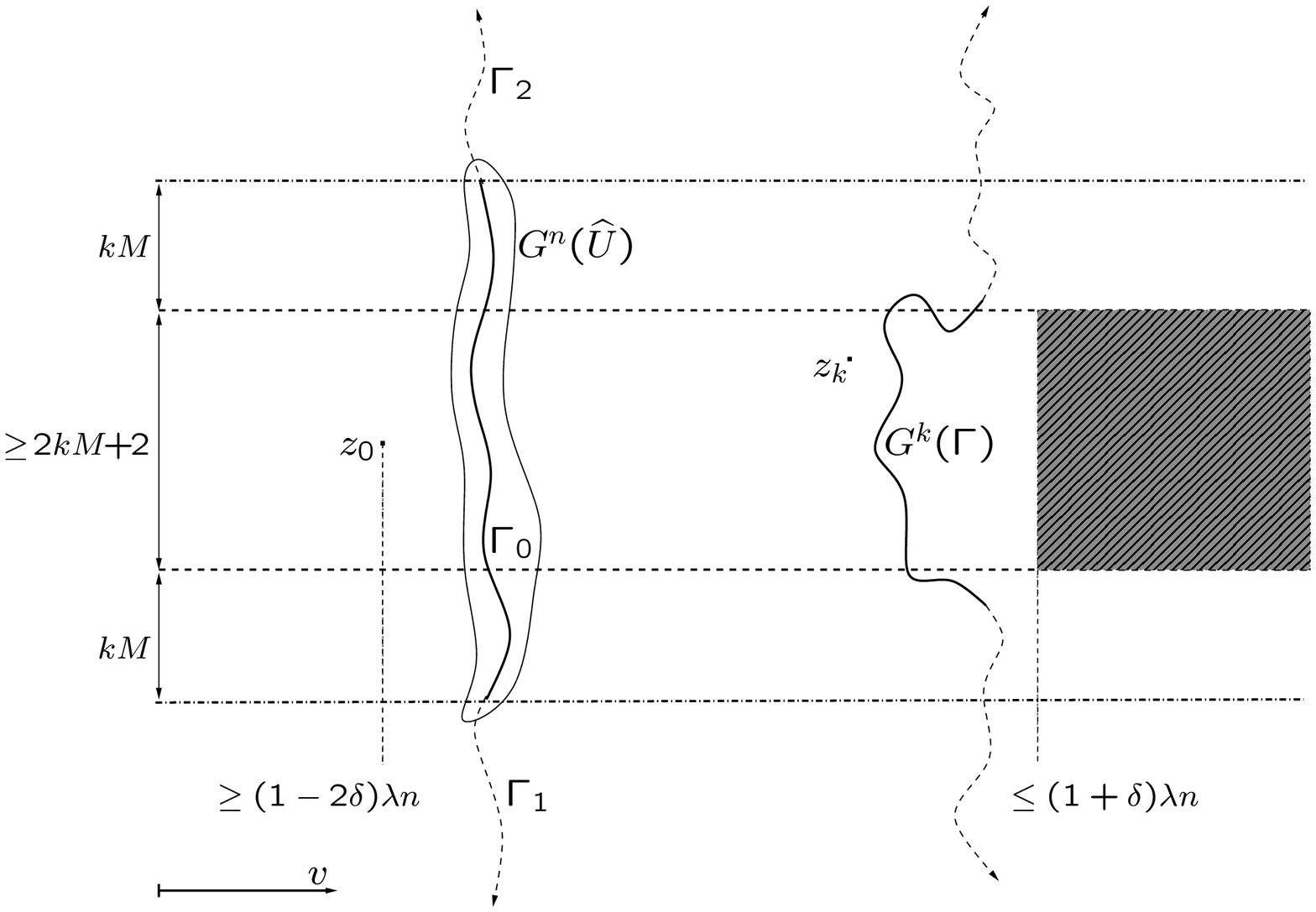, clip=,
      width=0.8\linewidth} \hspace{2eM} \end{center}
  \caption{\small Proof of Lemma~\ref{l.uniform-v-speed}: The slow movement of
    the curve $\Gamma$ impedes a faster movement of the point $z_0$ under
    iteration. \label{f.5} }
\end{figure}

Due to Corollary~\ref{c.linear-spreading} there exists a simple arc $\Gamma_0
\ssq G^n(\widehat{U})$ with endpoints $\zeta_1,\zeta_2\in G^n(\widehat{U})$ such
that $\vpbr{\zeta_2-\zeta_1} > cn$ and $\Gamma_0 \ssq
S_{v^\perp}[\vpbr{\zeta_1},\vpbr{\zeta_2}]$.  Let $\Gamma_1$ and $\Gamma_2$ be
properly embedded half-lines (proper images of $\R^+$) which join $\zeta_1$,
respectively $\zeta_2$, to infinity and satisfy
\begin{equation} \label{e.Gamma_i-unifspeed}
\Gamma_i \smin\{z_i\} \cap
S_{v^\perp}[\vpbr{\zeta_1},\vpbr{\zeta_2}] \  =  \ \emptyset \ .
\end{equation}
Let the properly embedded line $\Gamma = \Gamma_0 \cup \Gamma_1 \cup \Gamma_2$
be oriented in the direction from $\zeta_1$ to $\zeta_2$ and denote by $W$ the
connected component of $\R^2 \smin\Gamma$ to the left of $\Gamma$.  Then, since
$\Gamma_0\ssq G^n(\widehat{U})$ by assumption and due to (\ref{e.n_0-unifspeed})
and (\ref{e.Gamma_i-unifspeed}), $W$ contains the set
$S_v(-\infty,(1-\delta)\lambda n] \cap
S_{v^\perp}[\vpbr{\zeta_1},\vpbr{\zeta_2}]$.  Further, due to
(\ref{e.n_0-second}) and (\ref{e.k_0-n_0-ineq}) the set
\begin{equation} \label{e.W_0} W_0 \ = \ S_v[(1-2\delta)\lambda
  n,(1-\delta)\lambda n] \cap S_{v^\perp}[\vpbr{\zeta_1}+2kM,\vpbr{\zeta_2}-2kM]
\end{equation}
is a rectangle whose side-lengths are greater than 2. Hence $W_0$ contains a
fundamental domain of $\torus$ and by replacing $z_0$ with an integer translate
if necessary we may assume $z_0\in W_0$. This implies in particular that
$\vbr{z_0} \geq (1-2\delta)\lambda n$.  Now, there holds
\begin{equation}
  \label{e.z_k} z_k \ \in \ G^k(W) \cap S_{v^\perp}[\vpbr{\zeta_1} +kM,\vpbr{\zeta_2}-kM] \ . 
\end{equation}
However, due to (\ref{e.n_0-unifspeed}) and (\ref{e.Gamma_i-unifspeed}) we have 
\begin{equation}
G^k(W) \cap S_v((1+\delta)\lambda(n+k),\infty) \cap
S_{v^\perp}[\vpbr{\zeta_1} +kM,\vpbr{\zeta_2}-kM] \ = \ \emptyset \ ,
\end{equation}
such that $z_k \in S_v(-\infty,(1+\delta)\lambda(n+k)]$. Thus, using
(\ref{e.k_0-n_0-ineq}) and (\ref{e.delta-unifspeed}) we obtain
$\vbr{z_k-z_0} \leq (1+\delta)\lambda(n+k) - (1-2\delta)\lambda n <
(1+\eta)\lambda k$, in contradiction to (\ref{e.k_0-unifspeed}). \qed
\medskip

More or less along the same lines we obtain the following. 
\begin{lem}
  \label{l.rational-slope-segment} Suppose $\rho_{\widehat{U}}(G) \ssq \R \cdot
  v$ is a line segment of positive length.  Then $\rho(G) \ssq \R \cdot v$.
\end{lem}
\myproof Suppose for a contradiction that $\rho(G) \nsubseteq \R\cdot
v$ and assume withhout loss of generality that there exists
$\rho\in\textrm{Ex}(\rho(G))$ with $\sigma := \vpbr{\rho}>0$. Then as
above, there exists a point $z_0\in\R^2$ with $\nLim (z_n-z_0)/n =
\rho$, where $z_n=G^n(z_0)$. Hence, there exists an integer $N\geq 1$
such that
\begin{equation} \label{e.vertialspeed-in-rationalslope}
  \vpbr{z_n-z_0} \ > \ n\sigma/2 \quad \forall n \geq N \ .
\end{equation}
Let $av$ and $bv$ be the endpoints of $\rho_U(G)$, with $a<b$. Then by
definition of $\rho_{\widehat{U}}(G)$ in (\ref{e.rotation-subset}) and
due to the connectedness of $\bigcup_{j=0}^m G^j(\widehat{U})$, there
are infinitely many $m\in\N$ such that there exists a simple arc
$\Gamma_0 \ssq \bigcup_{j=0}^m G^j(\widehat{U})$ with endpoints
$\zeta_1$ and $\zeta_2$ such that
\begin{equation} \label{e.endpoint-distance-in-rationalslope}
  \vbr{\zeta_2-\zeta_1} \ \geq \ m(b-a)/3 \ . 
\end{equation}
Of course, we may assume $\Gamma_0 \ssq
S_v[\vbr{\zeta_1},\vbr{\zeta_2}]$.  Further, given any $\delta>0$
there exists $m_0\in\N$ such that for all $m\geq m_0$ we have
\begin{equation} \label{e.Gamma_0-orbit-in-rationalslope} G^n(\Gamma_0) \ \ssq \
  \bigcup_{j=0}^{n+m}G^j(\widehat U) \ \ssq \
  S_{v^\perp}[-\delta(m+n),\delta(m+n)] \quad \forall n\in\N \ .
\end{equation}
Now, fix $\delta > 0$ such that
\begin{equation}
  \label{e.delta-in-rationalslope}
  4\delta/(\sigma-2\delta) \ < \ (b-a)/12 \ 
\end{equation}
and choose $m\geq m_0$ such that there exists an integer $n_0\geq N$
with
\begin{equation}
  \label{e.i-in-rationalslope}
  n_0 \ \in \ \left[\frac{4\delta m +4}{\sigma - 2\delta},
    \frac{m(b-a)-6}{12M}\right] \ . 
\end{equation}
Let $\Gamma_1$ and $\Gamma_2$ be properly embedded half-lines such that
$\Gamma_i$ joins $\zeta_i$ to infinity and $\Gamma_i\smin\{\zeta_i\}$ is disjoint from
$S_v[\vbr{\zeta_1},\vbr{\zeta_2}]$. Let $\Gamma =
\Gamma_0\cup\Gamma_1\cup\Gamma_2$ be oriented in the direction from $\zeta_1$ to
$\zeta_2$. Denote by $W$ the connected component of $\R^2 \smin \Gamma$ to the
right of $\Gamma$. Then $W_0=S_v[\vbr{\zeta_1}+2n_0M,\vbr{\zeta_2}-2n_0M] \cap
S_{v^\perp}[-\delta m-2,-\delta m) \ssq W$ contains a fundamental domains of
$\T^2$. (Note that due to (\ref{e.endpoint-distance-in-rationalslope}) and
(\ref{e.i-in-rationalslope}) there holds $\vbr{\zeta_2}-\vbr{\zeta_1}\geq
4n_0M+2$.) Replacing $z_0$ by an integer translate if necessary, we may
therefore assume $z_0\in W_0$. It follows that
\begin{equation}
  z_{n_0} \ \in \ G^{n_0}(W) \cap S_v[\vbr{\zeta_1}+n_0M,\vbr{\zeta_2}-n_0M] \ .
\end{equation}
However, due to (\ref{e.Gamma_0-orbit-in-rationalslope}) and the choice of
$\Gamma$ the set $G^{n_0}(W)$ is disjoint from
$S_v[\vbr{\zeta_1}+n_0M,\vbr{\zeta_2}-n_0M] \cap
S_{v^\perp}(\delta(m+n_0),\infty)$, such that
\begin{equation}
  \label{e.z_i0-end}
  \vpbr{z_{n_0}-z_0} \ \leq \ \delta(2m+n_0)+2 \ . 
\end{equation}
However, from (\ref{e.i-in-rationalslope}) we obtain that $\delta(2m+n_0)+2 \leq
n_0\sigma/2$, such that (\ref{e.z_i0-end}) contradicts
(\ref{e.vertialspeed-in-rationalslope}).
\qed \medskip

\section{Proof of the main results}  \label{Proofs}

The following basic observation will be used in the proof of
Theorem~\ref{t.semilocal-rotsets}. Recall that $\varphi_n(z) =
(F^n(z)-z)/n$.

\begin{lem} \label{l.rotsets-convergence} Suppose $F$ is a lift of $f\in\homd$
  and $U\ssq \T^d$ is open, connected, bounded and recurrent.  Further, assume
  that $\rho_U(F) = {\cal S} \ssq \R \cdot v$ is a line segment with
  $0\notin{\cal S}$ and $v,v'$ are linearly independent.  Let $\widehat U$ be a
  connected component of $\pi^{-1}(U)$. Then there exists $p\in\N$ and $w \in
  \Z^2$ linearly independent of $v'$ such that
  \begin{equation} \label{e.pw-pairs}
   \left(F^p(\widehat U)- w\right)\cap \widehat U \ \neq \ 
  \emptyset \ .
\end{equation}
In particular, if $v$ is not the scalar multiple of an integer vector then there
exist infinitely many pairs $(p_i,w_i)$ with pairwise independent integer
vectors $w_i$ that satisfy (\ref{e.pw-pairs}). \end{lem} \proof Choose $\eps >0$
such that $B_{2\eps}(\rho_U(F)) = B_{2\eps}({\cal S})$ is disjoint from $\R
v'$. As $\widehat U$ is bounded, there exists $n_0\in\N$ such that
\[
\ntel\left(F^n(\widehat U)-\widehat U\right) \ \ssq \ B_\eps(\varphi_n(\widehat
U)) \quad \forall n\geq n_0 .
\]
(Here $A-B = \{ z-z' \mid z\in A,\ z'\in B \}$.)  Due to
Lemma~\ref{l.convergence} we may further assume, by increasing $n_0$ if
necessary, that
\[ 
\varphi_n(\widehat U) \ \ssq \ B_\eps(\rho_U(F)) \quad \forall n\geq n_0 \ .
\]
We thus obtain
\[
 \left(F^n(\widehat U)-\widehat U\right) \cap \R
  v' \ = \ \emptyset \quad \forall n\geq n_0 \ .
\]
Now, since $U$ is recurrent, there exists $p\geq n_0$ with $f^p(U)\cap U \neq
\emptyset$ and hence an integer vector $w$ with $\left(F^p(\widehat U)-w\right)
\cap \widehat U \neq \emptyset$. As $w$ belongs to $F^n(\widehat U)-\widehat U$
it must be linearly independent of $v'$.
 \qed\medskip

\noindent
{\bf Proof of Theorem~\ref{t.semilocal-rotsets}.} \ We start with (b) and
(c) and prove (a) at the end.  \alphlist \addtocounter{enumi}{1}
\item Let $f,F,U$ and ${\cal S} = \rho_U(F)$ be as in the statement of the
  theorem. Let $\widehat U$ be a connected component of $\pi^{-1}(U)$. We first
  assume that the line passing through ${\cal S}$ does not contain any rational
  points.  Since $U$ is non-wandering, there exists $p\in\N$ and $w\in\Z^2$ such
  that $\left(F^p(\widehat U)-w\right)\cap \widehat U \neq \emptyset$. Let $G =
  F^p-w$. Then $\rho_{\widehat U}(G) = p{\cal S}-w = L_{\lambda,v}[a,b]$ for
  some $\lambda\in\R,\ v\in\R^2_*$ and $a<b$. Further, the line $L_{\lambda,v}$
  contains no rational points either and therefore $\lambda \neq 0$. Hence
  $g=f^p,\ G$ and $\widehat U$ satisfy the assumptions of
  Lemma~\ref{l.uniform-v-speed} and we obtain $\rho(G)=L_{\lambda,v}[a,b]$ and
  thus $\rho(F) = (\rho(G)+w)/p ={\cal S}$.

  It remains to treat the case where the line passing through ${\cal S}$ has
  irrational slope and contains a single rational point. We choose $\widehat U$
  and $G$ as above and again have $\rho_{\widehat U}(G) = L_{\lambda,v}[a,b]$
  for some $\lambda\in\R,\ v\in\R^2_*$ and $a<b$. If $\lambda \neq 0$, then we
  can proceed exactly as before. However, this time we may have $\lambda = 0$
  since $L_{\lambda,v}$ may pass through 0. In this case
  Lemma~\ref{l.rotsets-convergence} yields the existence of a pair $\tilde
  p\in\N$ and $\tilde w\in\Z^2$ such that $\tilde w$ and $w$ are linearly
  independent and $\tilde G(\widehat U) \cap \widehat U \neq \emptyset$, where
  $\tilde G = F^{\tilde p}-\tilde w$. Then
\[
  \rho_{\widehat U}(\tilde G) \ = \ \tilde p\cdot \rho_{\widehat U}(F)-\tilde w \ = \
  \tilde p \cdot \left(\frac{\rho_{\widehat U}(G) + w}{p}\right) - \tilde w \ \ssq \ L_{0,v} + 
  \left(\frac{\tilde p}{p} \cdot w -\tilde w\right) 
\]
As $w$ and $\tilde w$ are linearly independent there holds $\frac{\tilde p}{p}
\cdot w -\tilde w \neq 0$. At the same time, this vector is not in
$L_{0,v}=\{v\}^\perp$ since the only rational vector contained in this line is
$0$. Therefore $\rho_{\widehat U}(\tilde G) = L_{\tilde \lambda, v}$ for some
$\tilde \lambda \neq 0$ and we can apply Lemma~\ref{l.uniform-v-speed} to
$\tilde g=f^{\tilde p},\tilde G$ and $\widehat U$ to obtain $\rho_{\widehat
  U}(F) = {\cal S}$ as before.
\item Suppose $\rho_U(F) = \{\rho\}$ with $\rho$ irrational. As above, we choose
  $p\in\N,\ w\in\Z^2$ and $G=F^p-w$ such that $G(\widehat U) \cap \widehat U
  \neq \emptyset$. Then $\rho_{\widehat U}(G) = \{p\rho-w\} =
  L_{1,p\rho-w}[0,0]$. By Lemma~\ref{l.extremalpoints} all extremal points of
  $\rho(G)$ are contained in $C_{p\rho-w}[0,0] = \R \cdot (p\rho-w)$ and hence
  $\rho(G) \ssq \R \cdot (p\rho-w)$. This implies $\rho(F) \ssq
  \R\cdot(\rho-w/p) + w/p =: A_1$.

  Due to Lemma~(\ref{e.pw-pairs}) we can repeat this argument with a second pair
  $\tilde p\in\N$ and $\tilde w \in \Z^2$, with $\tilde w$ linearly independent
  of $w$, and obtain $\rho(F) \ssq \R\cdot(\rho - \tilde w/\tilde p ) + \tilde
  w/\tilde p =: A_2$. It follows that $\rho(F)$ is contained in the intersection
  of the two lines $A_1$ and $A_2$, which is equal to $\{\rho\}$. (Note that as
  $\rho$ is irrational and the vectors $w/p$ and $\tilde w/\tilde p$ are both
  rational and linearly independent, the vectors $\rho-w/p$ and $\rho-\tilde
  w/\tilde p$ are linearly independent as well.)  \addtocounter{enumi}{-3}
\item Due to (a) and (b), it only remains to treat the two cases where
  $\rho_U(F)$ is either a line segment of positive length containing a rational
  vector or $\rho_U(F)$ is reduced to a single semi-rational vector. The second
  case is treated exactly as in (b), the only difference is that we cannot
  necessarily repeat the argument with a second $\tilde w$ to conclude that
  $\rho(F)$ is a singleton. 

  Hence, suppose that $\rho_U(F) \ssq \R \cdot v$ is a line segment of positive
  length and $\rho_U(F)$ contains a rational. By going over to a suitable
  iterate, we may assume w.l.o.g.\ that $0\in\rho_U(F)$. Again, there exists a
  pair $p\in\N$ and $w\in\Z^2$ such that $G(\widehat U) \cap \widehat U \neq 0$,
  where $G=F^p-w$. We have $\rho_{\widehat U}(G) \ssq L_{\lambda,v^\perp}$ with
  $\lambda := \vpbr{w}$. If $\lambda \neq 0$, then Lemma~\ref{l.uniform-v-speed}
  implies that that $\rho(G) = \rho_{\widehat U}(G)$ and hence $\rho(F) =
  \rho_U(F)$. If $\lambda = 0$, such that $\rho_{\widehat U}(G) \ssq \R\cdot v$,
  then Lemma~\ref{l.rational-slope-segment} implies $\rho(G) \ssq \R\cdot v$,
  such that we obtain $\rho(F) = \rho(G)/p \ssq \R\cdot v$ as well.  \qed
  \medskip \listend

\noindent
  {\bf Proof of Theorem~\ref{t.elliptic-islands}.}\  Suppose $U\ssq \torus$
  is open, bounded and connected. As $f$ in non-wandering, $U$ is also
  recurrent.  Suppose $\conv(\rho_U(F))$ has empty interior, that is,
  $\rho_U(F)$ is contained in a line. Then Theorem~\ref{t.semilocal-rotsets}
  implies that $\rho_U(F)$ is reduced to a single rational vector. By going over
  to a suitable iterate, we may assume $\rho_U(F) = \{0\}$. Let ${\cal D} := \{
  z\in\T^2\mid \exists \eps > 0 : \rho_{B_\eps(z)}(F) = 0\}$ and note that $U
  \ssq {\cal D}$. Let $D$ denote the connected component of ${\cal D}$ that
  contains $U$. Since ${\cal D}$ is $f$-invariant and $f$ is non-wandering, $D$
  is periodic with period $p$ for some $p\in\N$.

  We claim that $D$ contains no essential simple closed curve. In order to see
  this, suppose for a contradiction that the curve $\gamma \ssq D$ is essential
  with homotopy vector $v\in\Z^2_*$. Let $\Gamma$ be a connected component of
  $\pi^{-1}(\gamma)$. Then $\Gamma$ is a properly embedded line that remains in
  a bounded distance of $\R\cdot v$. Furthermore, due to compactness $\gamma\ssq
  D$ is covered by a finite number of open sets $U_i$ with $\rho_{U_i}(F) =
  \{0\}$, which implies $\rho_\gamma(F) = \{0\}$. However, this clearly
  contradicts the existence of rotation vectors $\rho$ with $\vbr{\rho} \neq 0$
  in $\rho(F)$. (Compare, for example, the proof of Lemma
  \ref{l.rational-slope-segment}.)

  In a similar way, we see that $D$ is simply-connected. If $\Gamma
  \ssq \pi^{-1}(D)$ is a closed Jordan curve, then by compactness we
  obtain $\rho_\Gamma(F) = \{0\}$. Consequently the Jordan domain
  $J(\Gamma)$ bounded by $J$ has rotation subset $\rho_{J(\Gamma)}(F)
  = \{0\}$ and therefore belongs to $\pi^{-1}(D)$ as well. Thus $D$ is
  the required $f^p$-invariant topological disk. Finally, as $D$ is
  homeomorphic to $\R^2$, the restriction $f^p_{|D}$ defines a plane
  homeomorphism. Since $f^p_{|D}$ is non-wandering, it must have
  a fixed point. (This follows, for example, from the Brouwer Plane
  Translation Theorem or from Franks Lemma \cite{franks:1988}).
  \qed\medskip

\noindent
  {\bf Proof of Proposition~\ref{p.chaotic}.} \ \alphlist \item Suppose
  $f\in\homtwo$ is non-wandering, $F$ is a lift and $z\in{\cal C}(f)$ is
  $\halb$-Lyapunov stable. Choose $\delta>0$ such that $f^n(B_\delta(z)) \ssq
  B_\halb(f^n(z)) \ \forall n\in\N$. Let $z_0\in\R^2$ be a lift of $z$ and $\wh
  U := B_\delta(z_0)$. Then $F^n(\wh U) \ssq B_\halb(F^n(z_0)) \ \forall
  n\in\N$, in particular $\diam\left(F^n(\wh U)\right) \leq 1$.

  Since $f$ is non-wandering $U$ is recurrent, such that there exist infinitely
  many pairs $(p,w) \in \N\times \Z^2$ with $\left(F^p(\wh U) - w\right) \cap
  \wh U \neq \emptyset$. The latter implies $F^p(\wh U) \ssq B_{1+\delta}(z_0 +
  w)$. Further, as $\left(F^{np}(\wh U) -w\right) \cap F^{(n-1)p}(\wh U) \neq
  \emptyset \ \forall n\in\N$, we obtain by induction that
  \[
  F^{np}(\wh U) \ \ssq \ B_{n+\delta}(z_0+nw) \quad \forall n\in\N \  . 
  \]
  However, this implies $\rho_{\wh U}(F) \ssq B_{\frac{1}{p}}(w)$, such that
  $\diam(\rho_{\wh U}(F)) \leq 1/p$. As $p$ can be chosen arbitrarily large we
  obtain $\diam(\rho_{\wh U}(F)) = 0$, in contradiction to $z\in{\cal C}(f)$.
\item Suppose $f$ is area-preserving, $U$ is a connected and bounded
  neighbourhood of $z\in{\cal C}(f)$ and $\wh U$ is a connected component of
  $\pi^{-1}(U)$. Birkhoffs Ergodic Theorem implies that Lebesgue-almost every
  point $z'\in\wh U$ has a rotation vector $\rho(F,z') =
  \nLim\left(F^n(z')-z'\right)/n$ (that is, the limit exists). If
  $\nLim\diam\left(\vbr{F^n(\wh U)}\right)/n=0$ this implies $\rho_{\wh U}(F)
  \ssq L_{\lambda,v}$ with $\lambda=\vbr{\rho(F,z')}$, contradicting $z\in {\cal
    C}(f)$.  \qed\listend

  \section{A parameter family of Misiurewicz-Ziemian type}
\label{MZ-family}

Examples of a toral homeomorphisms, homotopic to the identity, whose
rotation set has non-empty interior were introduced by Misiurewicz
and Ziemian \cite{misiurewicz/ziemian:1991} via lifts of the form
  \begin{equation}
    \label{e.mz-example}
    F(x,y) \ = \ (x+\psi_2(y+\psi_1(x)),y+\psi_1(x)) \ ,
  \end{equation}
  with continuous and 1-periodic functions $\psi_i:\R\to\R$, $i=1,2$.
  When $\psi_1(0)=\psi_2(0)=0$ and $\psi_1(\halb)=\psi_2(\halb)=1$, it
  can be easily checked that the points $(0,0),(\halb,0),(0,\halb)$
  and $(\halb,\halb)$ are fixed under the induced map $f\in\homtwo$
  and have rotation vectors $(0,0),(0,1),(1,0)$ and $(1,1)$,
  respectively.  Since the rotation set $\rho(F)$ is convex, it
  contains the square $[0,1]^2$.  More or less the same type of
  examples was proposed independently by Llibre and MacKay in
  \cite{llibre/mackay:1991}.

  In order to give a smooth example, we slightly modify this structure
  and let $\psi_1(x)=\alpha\sin(2\pi x)$ and $\psi_2(y)=\beta\sin(2\pi
  y)$ to obtain the parameter family
  \begin{equation}
    \label{e.mz-family}
    F_{\alpha,\beta}(x,y) \ = \ (x+\beta\sin(2\pi(y+\alpha\sin(2\pi x))),
    y+\alpha\sin(2\pi x)) \ . 
  \end{equation}
  We denote the toral diffeomorphisms induced by these lifts by
  $f_{\alpha,\beta}$. We note that $f_{\alpha,\beta}$ is area-preserving for
  all $\alpha,\beta\in\R$ and hence, in particular, non-wandering.

  For $F_*=F_{\halb,\halb}$, the points
  $(\viertel,0),(0,\viertel),(\dreiviertel,0)$ and $(0,\dreiviertel)$
  are $2$-periodic with rotation vectors
  $(0,\halb),(\halb,0),(0,-\halb)$ and $(-\halb,0)$. The rotation set
  $\rho(F_*)$ therefore contains the square $Q=\{(x,y)\mid |x|+|y|\leq
  \halb\}$ spanned by these vectors. Hence, the toral diffeomorphism
  $f_*=f_{\halb,\halb}$ induced by $F_*$ satisfies all the assumptions
  of Theorem~\ref{t.semilocal-rotsets}. Furthermore, this remains true
  for all $(\alpha,\beta)$ in a neighbourhood of $(\halb,\halb)$. The
  reason is that when the rotation set has non-empty interior, then it
  depends continuously on the toral homeomorphism in the ${\cal
    C}^0$-topology \cite{misiurewicz/ziemian:1991}.

\begin{figure}[h] 
\begin{center}
\epsfig{file=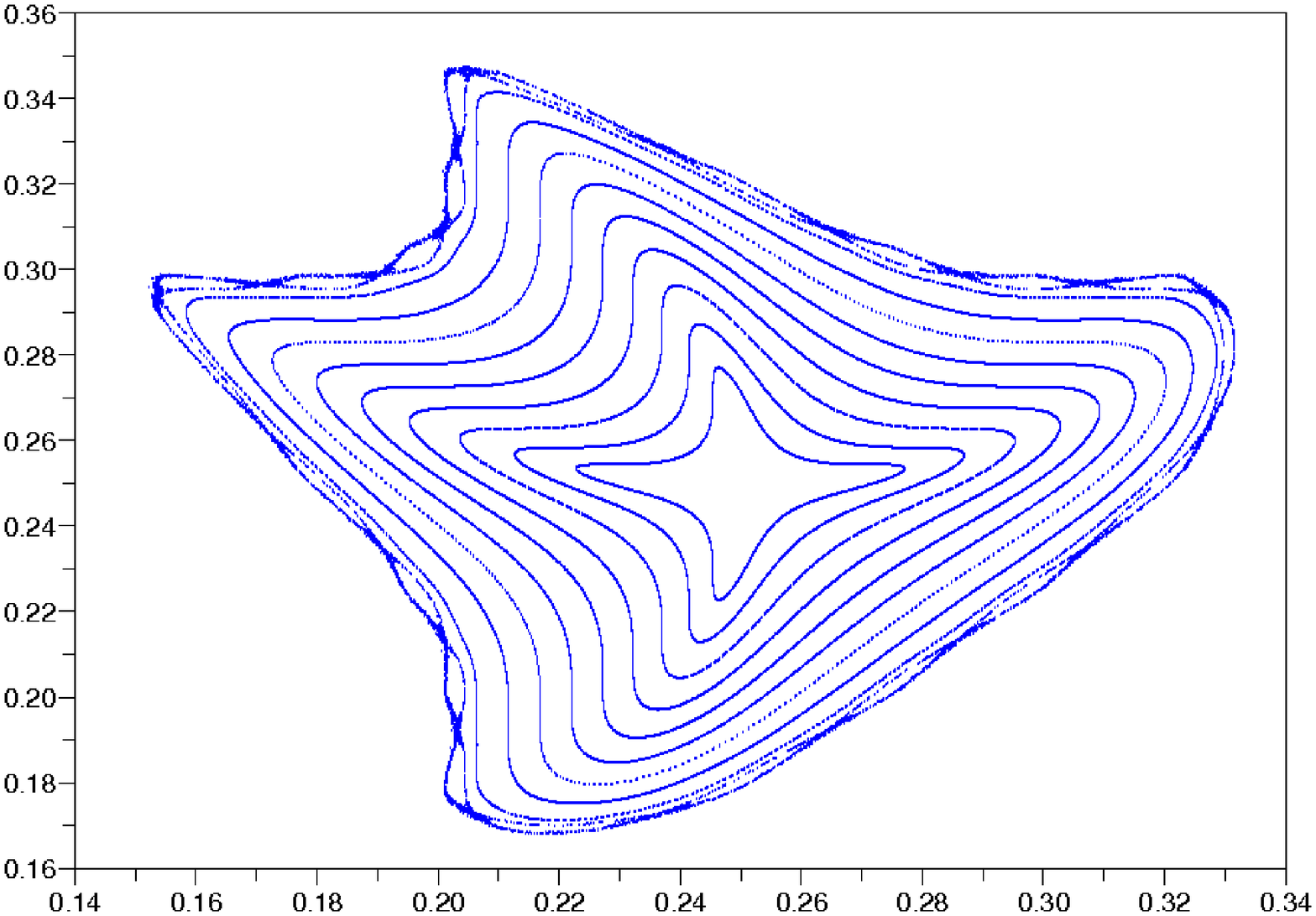, width=0.45\linewidth, height=0.45\linewidth}
\epsfig{file=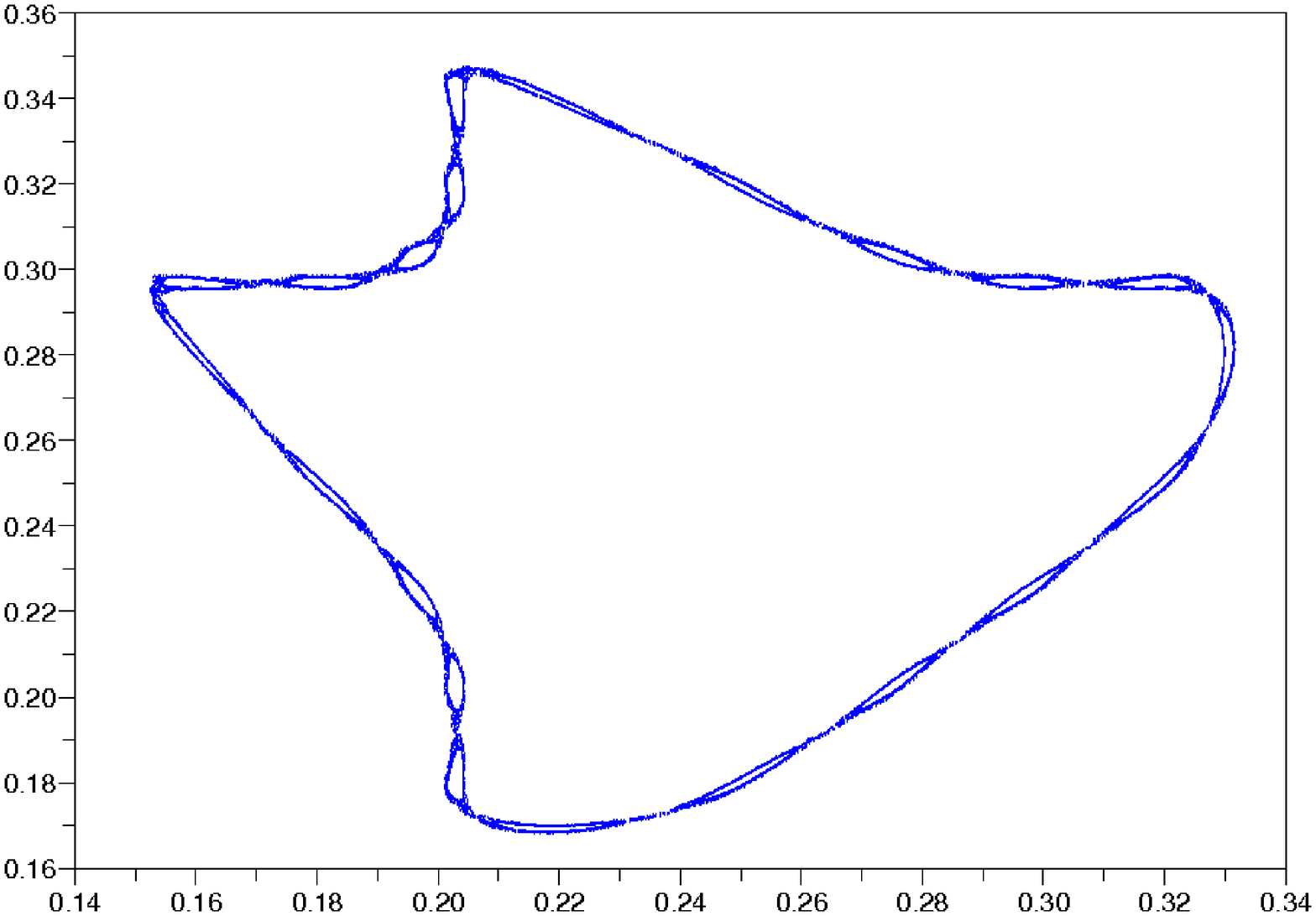, width=0.45\linewidth, height=0.45\linewidth}
\end{center} {\caption{ \small Enlargement of the elliptic island around
    $(\viertel,\viertel)$ in Figure~\ref{f.1}. Left: 20000 iterates of the
    starting points $(0.253+i\cdot 0.00455)\cdot(1,1),\ i=1\ld 10$. Right:
    100000 iterates of the starting point $0.298429 \cdot (1,1)$ in the
    instability region. ) \label{fig.7}}}
\end{figure}

  The points $(\frac{i}{4},\frac{j}{4})$ with $i,j=1,3$ are
  $2$-periodic for $f_*$ and have rotation vector $(0,0)$. The phase
  portrait of $f_*$ in Figure~\ref{f.1}\foot{Picture produced with
    SCILAB by computing 1000 iterates for each starting point in a
    40x40-grid on \torus\ and omitting the first 100 iterates for the
    plot.\label{foot.scilab}} clearly suggests that these points are
  surrounded by (star-shaped) elliptic islands. Since the differential
  $DF^2_*$ in these points is the identity matrix, it is difficult to
  establish the existence of elliptic islands in a rigorous way by the
  application of standard KAM-results%
  \foot{This would require an elliptic differential matrix with
    irrational rotation angle \cite{katok/hasselblatt:1997}.}, and we
  shall not further pursue this issue here. Yet, the pictures obtained
  by simulations show typical features of a KAM-type elliptic region,
  with invariant curves prevalent in the centre and at least one
  clearly visible instability zone towards the boundary (Figure~\ref{fig.7}).

Instead, we want to close by briefly discussing some of the symmetries present
in Figure~\ref{f.1}. The picture is invariant under the rotation with angle
$\pi/2$ around the point $(\halb,\halb)$, which is given by the map
$R(x,y)=(-y,x)$ on the torus. However, $f_*$ is not conjugate to itself by
$R$. The reason why this symmetry nevertheless appears is the fact that $R$
conjugates $f_*$ to its inverse, that is, $R^{-1}\circ f_*\circ R = f_*^{-1}$,
and for the visualisation of elliptic islands it does not make any difference if
$f_*$ is replaced by $f_*^{-1}$. This remains true for $f_{\alpha,\alpha}$ for all
$\alpha$, since
  \begin{equation}
    \label{e.sym}
    R^{-1}\circ f_{\alpha,\alpha}\circ R \ = \ f_{\alpha,\alpha}^{-1} \quad \forall \alpha\in\R\ . 
  \end{equation}
  Two symmetries conjugating $f_{\alpha,\beta}$ to itself, for all
  $\alpha,\beta\in\R$, are the rotation $S(x,y)=(-x,-y)$ with angle
  $\pi$ and the map $T(x,y) = (x+\halb,-y+\halb)$, which is the
  reflexion along the $x$-axis composed with the shift by
  $(\halb,\halb)$. This implies that if the points
  $(\frac{i}{4},\frac{j}{4})$ with $i,j=1,3$ are surrounded by
  elliptic islands, then these are all isometric to each other, but the
  isometries are orientation-reversing if $(i,j)$ and $(i',j')$ are
  such that $i=i'$ or $j=j'$. When $\alpha=\beta$ the additional
  symmetry given by $R$ implies that the elliptic islands are
  self-symmetric.  For the island surrounding $(\viertel,\viertel)$
  the symmetry axis is $L=\{(x,y)\mid x+y=\halb\}$, for the others it
  is the respective image of $L$ under $R^i,\ i=1,2,3$.

\begin{figure}[h] 
\begin{center} 
\epsfig{file=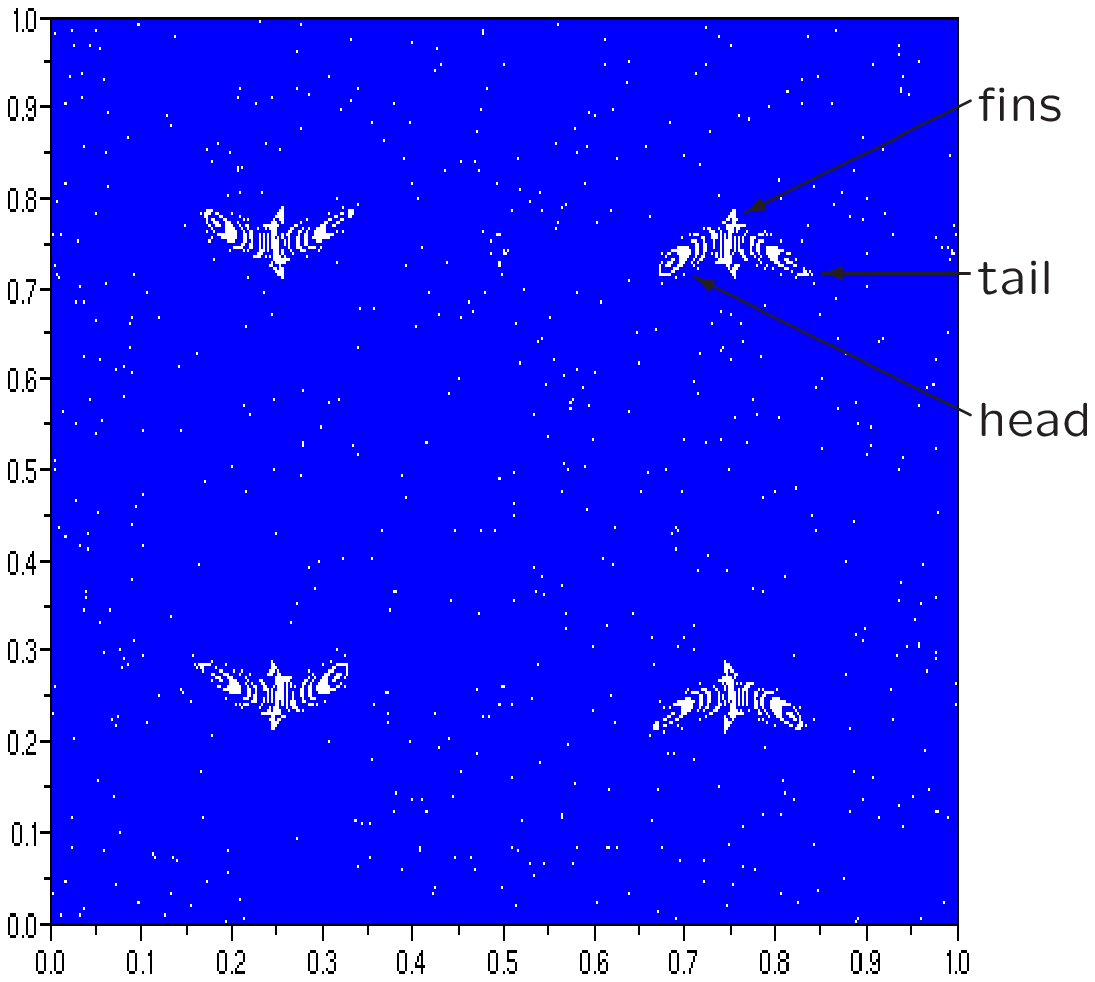, width=0.662\linewidth, height=0.58\linewidth}
\end{center} \vspace{-4ex} {\caption{ \small Elliptic sharks circling a chaotic ocean. (Phase
    portrait of $f_{\alpha,\beta}$ with $\alpha=0.5$ and
    $\beta=0.502$, with the same indications as in Footnote~\ref{foot.scilab}.) \label{f.6}}}
\end{figure}

  When $\alpha\neq \beta$ this self-symmetry of the elliptic islands
  breaks down since (\ref{e.sym}) does not hold anymore. With a slight
  adjustment of the parameters and some imagination in the Rorschach
  test below, this allows to return to a more aquatic environment
  (Figure~\ref{f.6}).

\vspace{-1ex}




\end{document}